\documentclass[12pt,a4paper, oneside, reqno]{amsart}

\usepackage[foot]{amsaddr} % authors’ affiliations just below the authors’ names
\usepackage{amsmath, nicefrac, amsthm, verbatim, amsfonts, amssymb, xcolor, bbm,mathtools}
\usepackage{graphics, xspace, enumerate}
\usepackage[a4paper,margin=2.5cm]{geometry}%this is for easy margins
\usepackage{todonotes}
\usepackage{dsfont}
\usepackage[bb=boondox]{mathalfa}

\usepackage{graphicx, amsmath}
\usepackage[colorlinks=true,citecolor=blue,urlcolor=green,linkcolor=blue,
bookmarksopen=true,unicode=true,pdffitwindow=false, pdfdisplaydoctitle]{hyperref}
\usepackage[english]{babel}
\usepackage[languagenames,fixlanguage]{babelbib}

\usepackage{seqsplit,cleveref}

\hypersetup{pdfauthor={}}

\theoremstyle{plain}
\newtheorem{theorem}{Theorem}[section]
\newtheorem{corollary}[theorem]{Corollary}
\newtheorem{lemma}[theorem]{Lemma}

\theoremstyle{definition}
\newtheorem{remark}[theorem]{Remark}

\newtheorem{assumption}{A\!\!}

\newcommand {\Prob} {\ensuremath{\mathbb{P}}}
\newcommand {\R} {\ensuremath{\mathbb{R}}}
\newcommand {\ZZ} {\ensuremath{\mathbb{Z}}}
\newcommand {\N} {\ensuremath{\mathbb{N}}}

\newcommand{\D}{\mathrm{d}}

\newcommand{\Var}{\operatorname{Var}}

\newcommand{\bbN}{\mathbb{N}}
\newcommand{\bbZ}{\mathbb{Z}}
\newcommand{\bbR}{\mathbb{R}}

\newcommand{\bbP}{\mathbb{P}}
\newcommand{\bbE}{\mathbb{E}}

\newcommand{\Ind}{\mathds{1}}
\newcommand{\calC}{\mathcal{C}}
\newcommand{\barcalC}{\overline{\mathcal{C}}}

\newcommand{\calF}{\mathcal{F}}

\newcommand{\calR}{\mathcal{R}}

\newcommand{\Capacity}{\mathrm{Cap}}

\newcommand{\aps}[1]{\vert #1 \vert}
\newcommand{\APS}[1]{\left\vert #1 \right\vert}
\newcommand{\norm}[1]{\lVert #1 \rVert}
\newcommand{\floor}[1]{\lfloor #1 \rfloor}

\numberwithin{equation}{section}

\title[Invariance principle for the range of stable random walks]{Invariance principle \\for the capacity and the cardinality\\ of the range of stable random walks}

\author[W.\ Cygan]{Wojciech Cygan $^{1,2}$}
\address{$^{1}$University of Wroc\l{}aw,
		Faculty of Mathematics and Computer Science\\
		Institute of Mathematics,
		pl.\ Grunwaldzki 2/4, 50--384 Wroc\l{}aw, Poland}
\address{$^{2}$Technische Universit\"{a}t Dresden,
		Faculty of Mathematics\\
		Institute of Mathematical Stochastics,
		Zellescher Weg 25, 01069 Dresden, Germany}
\email{wojciech.cygan@uwr.edu.pl, wojciech.cygan@tu-dresden.de}

\author[N.\ Sandri\'{c}]{Nikola Sandri\'{c} $^{3}$}
\address{$^{3}$Department of Mathematics\\University of Zagreb\\ Zagreb\\Croatia}
\email{nsandric@math.hr}

\author[S.\ \v{S}ebek]{Stjepan\ \v{S}ebek $^{4}$}
\address{$^{4}$Department of Applied Mathematics\\
	Faculty of Electrical Engineering and Computing\\
	University of Zagreb\\ 
 Zagreb\\ 
	Croatia}
\email{stjepan.sebek@fer.hr}


\begin{document}
\allowdisplaybreaks[4]

\begin{abstract}
We prove an almost sure invariance principle for the capacity and the cardinality of the range of a class of $\alpha$-stable random walks on the integer lattice $\ZZ^d$ with $d/\alpha > 5/2$, and $d/\alpha >3/2$, respectively. As a direct consequence, we conclude Khintchine's and Chung's laws of the iterated logarithm for both processes.
\end{abstract}

\maketitle

\section{Introduction}
Let $\{X_i\}_{i\ge1}$ be a sequence of independent and identically distributed random variables with values in $\bbZ^d$, defined on a given probability space $(\Omega,\calF,\bbP)$. We consider a random walk $S_n=X_1+\cdots  +X_n$, $n\ge1$, with  $S_0=0$. 
The range  $\{\calR_n\}_{n\ge0}$ of the random walk  $\{S_n\}_{n\ge0}$ is defined as the random set
\begin{equation*}
\calR_n = \{S_0,\dots,S_n\},\qquad n\ge 0,
\end{equation*}
and the cardinality of the range is denoted by $|\calR_n|$, $n\ge0$.

The capacity of a set $A\subseteq \mathbb{Z}^d$
(with respect to a transient random walk $\{S_n\}_{n\ge0}$) is defined as 
\begin{equation*}
{\rm Cap}(A)=\sum_{x\in A}\bbP(\tau^x_A=\infty),
\end{equation*}
where $\tau^x_A$ is the first return time of $\{S_n+x\}_{n\ge0}$ to the set $A$, that is, 
\begin{equation*}
\tau^x_A=\inf\{n\ge1:S_n+x\in A\}.
\end{equation*}
We use notation $\tau^x_y$ if $A=\{y\}$, and we abandon the upper index if $x=0$.

The main aim of this article is to investigate the growth of the capacity and the cardinality of the range process. More precisely, we consider a class of symmetric stable random walks  (see assumptions \textbf{A}\ref{A1} and \textbf{A}\ref{A2} below) and obtain an almost sure invariance principle  which asserts that the centered stochastic process $\{\Capacity (\calR_n) - \mathbb{E} [\Capacity (\calR_n)]\}_{n\geq 0}$ can be approximated (up to a constant) by a path of a standard Brownian motion almost surely (see Theorem \ref{TM:1} below). 
As a corollary, we obtain Khintchine's and Chung's laws of the iterated logarithm.
Our approach is based upon two main ingredients. The first is a decomposition of the capacity of the range which allows us to represent it as a sum of finitely many independent random variables plus an error term. The second is the Skorohod embedding theorem which we apply to replace  the sequence of independent random variables that appears in the capacity decomposition with a Brownian path  sampled at certain random instances.

A slick adjustment of our method  allows us to obtain analogous results for the process $\{|\calR_n|\}_{n\ge0}$ (see Theorem \ref{TM:2} below). For this case we utilize a decomposition of the range
into a sum of finitely many independent parts plus an error term which can be treated as the number of intersection points of two independent copies of the underlying random walk.

Problems related to the range of random walks constitute a rich area of research in modern probability theory. The first result in this direction is due to
Dvoretzky and Erd\"{o}s \cite{Dvoretzky}  where they obtained  strong law of large numbers for $\{|\calR_n|\}_{n\ge0}$ of a simple random walk in  $d\ge2$. Later on, Spitzer  \cite{Spitzer} extended this result to all random walks  in $d\ge1$.  A central limit theorem (CLT) for $\{|\calR_n|\}_{n\ge0}$ (with a normal law in the limit) was first  obtained by Jain and Orey \cite{Jain-Orey} for strongly transient random walks, and Jain and Pruit \cite{Jain_Pruitt}
extended this result to all random walks in $d\ge3$. 
Le Gall \cite{LeGall-French} proved a version of a CLT for $\{|\calR_n|\}_{n\ge0}$ of all two-dimensional random walks with zero mean and finite second moment with a non-normal law in the limit. 
Le Gall and Rosen \cite{LeGall-Rosen} established  strong law of large numbers and  CLT for $\{|\calR_n|\}_{n\ge0}$ of a class of $\alpha$-stable random walks in $\mathbb{Z}^d$. The law of the random variable in the limit depends on the value of the ratio $d/\alpha$ in this case. We also refer to \cite{Kang-Wee} for a distributional invariance principle for the range of one dimensional recurrent stable random walks.

The capacity of the range has attracted much attention in the literature as well. 
To understand the motivation for the study of the capacity of the range and its links to the theory of intersection of paths of random walks we refer the reader to \cite{Asselah_deviation_Cap} and  interesting references therein. The first result concerning the limiting behavior of the process $\{{\rm Cap}(\calR_n )\}_{n\ge0}$ is due to Jain and Orey \cite{Jain-Orey} where they obtained  strong law of large numbers for all transient random walks. CLT was recently proved in \cite{Asselah_Zd} for a simple random walk in $\ZZ^d$ with $d\geq 6$. The case of a simple random walk in dimensions $d=4$ and $d=5$ was studied in \cite{Asselah_Z4} and \cite{Schapira} respectively, see also \cite{Chang} for $d=3$. In \cite{CSS19} authors of the present article established CLT for the capacity of the range of a class of stable random walks which possess one-step loops, see also \cite{CSS19_2} for a functional CLT for such random walks. We refer to  \cite{Bai_2}, \cite{Bai_3} and  \cite{Bai_1}  for interesting results on capacity of the range of branching and tree-indexed random walks. Recently, a strong law of large numbers and CLT were obtained in \cite{RSS} for the capacity of the range of symmetric random walks on finitely generated groups of order $d$, for $d\geq 6$.

The theory of growth of paths of the process $\{{\rm Cap}(\calR_n )\}_{n\ge0}$ is still in its infancy. In \cite{Asselah_deviation_Cap} the authors proved downward large and moderate deviation estimates and an upward large deviation principle for a symmetric simple random walk in dimensions $d\geq 5$. Almost sure invariance principles and laws of the iterated logarithm for the capacity of the range of symmetric simple random walk in $\bbZ^d$ were recently obtained in \cite{Dembo-Okada} for dimensions $d\geq 3$. 
In this article, we establish analogous results for symmetric stable random walks which admit one-step loops (see assumption \textbf{A}\ref{A2} below), under the assumption that $d/\alpha> 5/2$, where $\alpha \in (0,2]$ is the index of stability.
We remark that the limit behavior of the capacity (and cardinality) of the range of stable random walks depends on the value of the ratio $d/\alpha$. 
One can observe an analogy between the limit behavior of the cardinality of the range for $d/\alpha >3/2$, and the capacity of the range for $d/\alpha>5/2$
(as already commented in \cite{CSS19}). If $d/\alpha \le 5/2$, we then  expect that the capacity of the range behaves differently than displayed in the present article as this is the case for the cardinality of the range when $d/\alpha\leq 3/2$.
We conjecture that if $d/\alpha=5/2$ then the Brownian path which appears in our almost sure invariance principle has to be evaluated at certain times $\phi (n)$ for a specific function $\phi(t)$ which is determined by the variance of the capacity of the range, see \cite{Bass_Kumagai}. Similarly, for $d/\alpha<5/2$ we suspect that the Brownian motion has to be replaced with another stochastic process, see  \cite{Bass_Rosen} where such a process is given by renormalized intersection local times of the Brownian motion in the case of the cardinality of the range of planar random walks.

As  it was mentioned before, our approach which we apply for the capacity of the range can be slightly changed and then transferred to obtain an almost sure invariance principle for $\{|\calR_n|\}_{n\ge0}$ when $d/\alpha>3 /2$.
The study of the growth of the process $\{|\calR_n|\}_{n\ge0}$ was initiated with Khintchine's law of the  iterated logarithm  obtained by Jain and Pruitt \cite{Jain_Pruitt_LIL} for aperiodic random walks satisfying $\mathbb{P}(\tau_0=\infty)<1$  and such that they are either strongly transient or lie in $\ZZ^d$ with $d\geq 4$.  We note that if $\mathbb{P}(\tau_0=\infty)=1$ then $|\calR_n| =n+1$ a.s.
Under similar assumptions Hamana \cite{Hamana} proved an almost sure invariance principle.
Bass and Kumagai \cite{Bass_Kumagai} obtained an  almost sure invariance principle,   and Khintchine's and Chung's laws of the  iterated logarithm for a class of random walks which have finite moments of order $2+\delta$, for $\delta>0$, in dimensions $d=2$ and $d=3$, see also \cite{Bass_Chen_Rosen} for a significant extension to random walks with finite second moments. The one-dimensional case was studied by Chen in \cite{Chen}. 
Our contribution to this topic is the study of random walks satisfying \textbf{A}\ref{A1} which clearly do not have to posses finite second moments nor finite supports. It means that by an appropriate choice of a small $\alpha$ we can achieve a rich class of random walks for which the almost sure invariance principle holds in all dimensions.
To efficiently construct examples of random walks satisfying \textbf{A}\ref{A1} one can employ the method of discrete subordination which was developed in \cite{BSC}. We refer to \cite{BCT} and \cite{Mimica}   for a condition under which a subordinate random walk belongs to the domain of attraction of a stable law. Moreover, random walks constructed according to this procedure fulfill \textbf{A}\ref{A2}.

\subsection*{Assumptions and main results}
In the course of study we confine our attention to aperiodic random walks only. The random walk $\{S_n\}_{n\ge0}$ is   aperiodic if the smallest additive subgroup generated by the set $\{x\in\mathbb{Z}^d:\, \Prob(S_1=x)>0\}$ is equal to $\ZZ^d$. This assumption  is not restrictive, for if $\{S_n\}_{n\ge0}$ is not aperiodic, we could then pose the problem (and prove the same theorems) on a smaller subgroup of $\ZZ^d$, see \cite[pp 20]{Spitzer}. 

We obtain our results for a class of symmetric $\alpha$-stable random walks in $\ZZ^d$, that is, we assume the following condition.

\begin{assumption}\label{A1}
The random walk $\{S_n\}_{n\ge0}$ is symmetric and it belongs to the domain of  attraction of a  non-degenerate $\alpha$-stable random law with $0<\alpha\le2$, meaning that there exists a regularly varying function $b(x)$ with index $1/\alpha$ such that 
	$$\frac{S_n}{b(n)}\xrightarrow[n\nearrow\infty]{\text{(d)}}X_\alpha,$$ 
	where $X_\alpha$ is an  $\alpha$-stable random variable on $\R^d$ and $\xrightarrow[]{\text{(d)}}$ stands for the convergence in distribution.
\end{assumption}

To perform a necessary asymptotic analysis of the variance of $\Capacity (\calR_n)$ we need an additional assumption.

\begin{assumption}\label{A2}
The random walk $\{S_n\}_{n\ge0}$ admits one-step loops, that is, $\Prob(X_1=0)>0$.
\end{assumption}

\smallskip

We now state the main results of the article.\footnote{We use standard  Landau notation: for $f:\N\to\R$ and $g:\N\to(0,\infty)$ we write $f(n)=O(g(n))$ if $\limsup_{n\to \infty } |f(n)|/g(n) <\infty$. Similarly, $f(n) = o(g(n))$ if $\lim_{n\to \infty}f(n)/g(n) = 0$. We also write $o_n(1)$ for a sequence converging to zero as $n$ goes to infinity.}

\begin{theorem}\label{TM:1}
	Assume \textbf{A}\ref{A1}, \textbf{A}\ref{A2} and $d/\alpha>5/2$. Then, there
	exists a standard Brownian motion 	$\{B_t\}_{t\ge0}$ (defined on the same, possibly enlarged, probability space as $\{S_n\}_{n\ge0}$) and 
	a constant $\sigma>0$, such that 
	\begin{equation}\label{TME:1.2}
	\sigma^{-1} \left( {\rm Cap}(\calR_n )-\mathbb{E}[{\rm Cap}(\calR_n )] \right) -B_n 
	=
	O\bigl(\sqrt{n}\bigr)
 \quad \text{a.s.}
	\end{equation}
	
\end{theorem}

With the aid of classical growth results for paths of standard Brownian motion (see e.g.\ \cite[Chapter 11]{schilling-book}), we conclude  Khintchine's and Chung's laws of the iterated logarithm.

\begin{corollary}
	Under the assumptions of \Cref{TM:1}, it holds $\mathbb{P}$-a.s.
	\begin{align*}
	&\liminf_{n\nearrow\infty}\frac{{\rm Cap}(\calR_n )-\mathbb{E}[{\rm Cap}(\calR_n )]}{\sqrt{n\log\log n}}=-\sqrt{2}\sigma,\\
	&\limsup_{n\nearrow\infty}\frac{{\rm Cap}(\calR_n )-\mathbb{E}[{\rm Cap}(\calR_n )]}{\sqrt{n\log\log n}}=\sqrt{2}\sigma,\\
	&\liminf_{n\nearrow\infty}\frac{\max_{0\le m\le n}\bigl|{\rm Cap}(\calR_m )-\mathbb{E}[{\rm Cap}(\calR_m )]\bigr|}{\sqrt{n/\log\log n}}=\frac{\pi}{8}\sigma.
	\end{align*}
\end{corollary}

Let us remark that assumption $d/\alpha >5/2$ implies that the random walk $\{S_n\}_{n\ge0}$ is strongly transient. We recall that a transient random walk $\{S_n\}_{n\ge0}$ is called strongly transient if $\sum_{n\ge1}n\,\Prob(S_n=0)<\infty$; otherwise it is called weakly transient.   
%It is known that every transient random walk is either strongly or weakly transient (see \cite{Sato}). 
According to \cite[Theorem 3.4]{Sato} and  \cite[Theorem 7]{Takeuchi}, each random walk  satisfying \textbf{A}\ref{A1} is transient if $d/\alpha>1$ and strongly transient if $d/\alpha>2$.  
We remark that  (strong) transience assumption is quite natural in the present context as it  ensures  that the range process $\{\calR_n\}_{n\ge0}$ grows fast enough and together with assumption \textbf{A}\ref{A2} it enables us to control the  variance of $\Capacity (\calR_n)$. More precisely, it was proved in \cite[Lemmas 4.3 and 5.3]{CSS19} that under \textbf{A}\ref{A1}-\textbf{A}\ref{A2} and for $d/\alpha >5/2$ there is a positive constant $\sigma$ such that
\begin{align}\label{Var-asymp}
\lim_{n\to \infty}\frac{\Var (\mathrm{Cap}(\calR_n))}{n}=\sigma >0.
\end{align}
The constant $\sigma$ from \eqref{Var-asymp} coincides with the constant appearing in Theorem \ref{TM:1}.
In this article, we establish the second order term in \eqref{Var-asymp}, see Lemma \ref{LM:3.2}.

Our strategy to prove \Cref{TM:1} is based upon a capacity decomposition recently established in \cite{Dembo-Okada} which bases upon the approach of \cite{Jain_Pruitt_LIL} and \cite{Bass_Kumagai}. This allows us to represent the random variable $\Capacity (\calR_n)$ as a sum of finitely many independent random variables plus an error term. To show that the error term is negligible we apply estimates for its moments from \cite{CSS19}. We then employ the Skorohod embedding theorem to approximate the sum of independent random variables  by a path of Brownian motion evaluated at a random time which is given by a sum of specific stopping times. 
To get rid of randomness coming from this sequence we make use of the second order asymptotics of $\Var (\Capacity (\calR_n))$. To find the second order term in \eqref{Var-asymp}, we apply another capacity decomposition extracted from \cite{Asselah_Zd}, see Section \ref{Appendix}.

Performing an analogous approach for the cardinality of the range, we establish the following result. 
\begin{theorem}\label{TM:2}
	Assume \textbf{A}\ref{A1} and  $d/\alpha>3/2$. 
	Then, there
	exists a standard Brownian motion 	$\{B_t\}_{t\ge0}$ (defined on the same, possibly enlarged, probability space  as $\{S_n\}_{n\ge0}$) and a constant $\tilde{\sigma}>0$, such that 
	\begin{equation*}
	\tilde{\sigma}^{-1} \left( |\calR_n |-\mathbb{E}[|\calR_n |]\right) -B_n
	=
	O\bigl(\sqrt{n}\bigr)\ \quad \text{a.s.}
	\end{equation*}
\end{theorem}
Similarly as before, we conclude the corresponding laws of the iterated logarithm.
\begin{corollary}
	Under the assumptions of \Cref{TM:2} it holds $\mathbb{P}$-a.s.
	\begin{align*}
	&\liminf_{n\nearrow\infty}\frac{|\calR_n |-\mathbb{E}[|\calR_n |]}{\sqrt{n\log\log n}}=-\sqrt{2}\tilde{\sigma},\\
	&\limsup_{n\nearrow\infty}\frac{|\calR_n |-\mathbb{E}[|\calR_n |]}{\sqrt{n\log\log n}}=\sqrt{2}\tilde{\sigma},\\
	&\liminf_{n\nearrow\infty}\frac{\max_{0\le m\le n}\bigl||\calR_m |-\mathbb{E}[|\calR_m |]\bigr|}{\sqrt{n/\log\log n}}=\frac{\pi}{8}\tilde{\sigma} .
	\end{align*}
\end{corollary}

To prove  \Cref{TM:2} we utilize a method of splitting the range into a sum of finitely many independent parts plus an error term  which can be treated as the number of intersection points of two independent copies of the underlying random walk. To deal with the error term in this case we apply estimates of moments of the number of intersection points which we extract from \cite{LeGall-Rosen}, see  \Cref{Sec:cardinality}. We then again employ the Skorohod embedding theorem and approximate the sum of independent random variables  by a path of Brownian motion evaluated at a random time. 
In this case we apply the second order asymptotics obtained in Lemma \ref{lm:var_asymp-card} which takes the form  $\Var (|\calR_n|)=\tilde{\sigma}n+O(n^{1/2}H(n))$, for a function $H(n)$ defined in \eqref{Function_H}. The constant $\tilde{\sigma}$ coincides with the constant appearing in Theorem \ref{TM:2}, and according to \cite[Theorem 4.4]{LeGall-Rosen} it is strictly positive. 

Let us briefly display the structure of the remaining part of the article. In Section \ref{Sec:capacity}, we focus on the proof of \Cref{TM:1} which is accompanied by a series of auxiliary results, see \Cref{Lemma-supV}, \Cref{lemma_difference} and \Cref{lemma_C-diff-max-estimate}.
In Section \ref{Sec:cardinality}, we deal with the cardinality process and we present the proof of \Cref{TM:2}. In the closing appendix (see Section \ref{Appendix}), we find the second order term for the variance asymptotics given at \eqref{Var-asymp} and we provide a corresponding result for the cardinality process, see \Cref{LM:3.2} and \Cref{lm:var_asymp-card}.

\section{The Capacity Process}\label{Sec:capacity}

This section is devoted to the proof of \Cref{TM:1}.
We combine the approach of \cite[Theorem 2.1]{Bass_Kumagai} together with that of \cite[Section 5]{Dembo-Okada} but let us emphasize that it requires a subtle analysis and numerous adjustments in the present setting.
For any $r,s,t \in \bbZ$, $0\leq r\leq  s\leq  t$ we use the notation $\calR_r = \{S_0,\ldots ,S_r\}$, $\calR_\infty = \{S_0, S_1, \ldots \}$,  $\calR(r,s] = \{S_{r+1},\ldots ,S_s\}$ (with the convention $\calR(r,r]=\emptyset$), $\calR (r, \infty ) = \{S_{r+1},S_{r+2},\ldots \}$, and respectively 
\begin{align*}
\calC_r = \mathrm{Cap}(\calR_r),\quad
 \calC(r,s] = \mathrm{Cap} (\calR (r,s]),\quad \text{and}\quad V_{r,s,t}=\calC(r,s]+\calC(s,t]-\calC(r,t].
\end{align*}
By translation invariance of the capacity, $\calC(r,s]$ has the same law as $\calC(0,s-r]$.  By independence of the increments of the random walk combined with  translation invariance of the capacity, the random variables $\calC(r,s]$ and $\calC(s,t]$ are independent. For any increasing sequence $\{n_k\}_{k\ge0}\subset \bbZ$ with $n_0=0$ we easily derive the following decomposition (see also \cite[Eq.\ (4.2)]{Dembo-Okada})
\begin{align}\label{Cap-decomp-Dembo-Okada}
\calC (0, n_k]
=
\sum_{j=1}^k U_j - \Delta_{n_k,k},
\end{align}
where $U_j = \calC(n_{j-1},n_j]$, $j=1,\dots,k$, are independent random variables and
\begin{align*}
 \Delta_{n_k,k} = \sum_{j=1}^{k-1} V_{n_{j-1},n_j,n_k}
 =
 \sum_{j=1}^{k-1} V_{0,n_j,n_{j+1}}.
\end{align*}

Recall that the function $b(x)$ appearing in assumption \textbf{A}\ref{A1} is necessarily of the form
\begin{align}\label{b_function}
b(x) = x^{1/\alpha} \ell(x),\qquad x>0,
\end{align}
where $\ell(x)$ is a slowly varying function, see \cite{BGT_book} for a detailed discussion on regularly varying functions. 
 In the sequel we will frequently use the following function
	\begin{equation}\label{h_d_function}
		h(n) = 
		\begin{cases}
			1, &\quad d/\alpha > 3 ;\\
			\sum_{k=1}^n k^{-1}(\ell (k))^{-d}, & \quad d/\alpha  = 3;\\
			n^{3-d/\alpha}(\ell(n))^{-d}, & \quad 2<d/\alpha <3.
		\end{cases}
	\end{equation} 
	In particular, $h(n)$ is slowly varying and increasing for $d/\alpha \geq 3$ (see \cite[Lemma 2.2]{LeGall-Rosen}), and regularly varying of index $3-d/\alpha$ for $d/\alpha \in (2,3)$. Since in the latter case the index of regular variation is positive, we can assume that $h(n) $ is non-decreasing, see \cite[Theorem 1.5.3]{BGT_book}.
	
By $G(x,y) $ we denote the Green function of the random walk $\{S_n\}_{n\geq 0}$, that is, 
\begin{align*}
G(x,y) \,=\, \sum_{n=0}^\infty \bbP (S_n=y-x)\,, \qquad x, y \in \bbZ^d\,.
\end{align*}
Also, for $A,B\subseteq\ZZ^d$ we use notation
\begin{align*}
G(A, B) \,=\, \sum_{x\in A}\sum_{y\in B} G(x,y)\,.
\end{align*}
We evidently have $G(A_1,B)\leq G(A_2,B)$ for $A_1\subseteq A_2$.

We start with a series of crucial lemmas.
\begin{lemma}\label{Lemma-supV}
Assume \textbf{A}\ref{A1} and $d/\alpha >2$.
For any $p\geq 1$ there exists a constant $c=c(d,\alpha ,p)>0$ such that for any integers $0\leq r < s<t$ it holds
 \footnote{For a random variable $Y$ we write $\|Y\|_p=(\mathbb{E}[|Y|^p])^{1/p}$.} 
\begin{align*}
\Vert \sup_{t> s}V_{r,s,t} \Vert_{p} \leq c\, h(s-r)\qquad \mathrm{and}\qquad
\Vert \sup_{0\leq r< s}V_{r,s,t} \Vert_p
\leq c\, h(t-s).
 \end{align*}
\end{lemma}

\begin{proof}
%By translation invariance of the capacity, we can assume that $r=0$. 
By \cite[Proposition 1.2]{Asselah_Zd}, 
\begin{align}\label{Assel-Cap-ineq}
V_{r,s,t}\leq 2 G(\calR (r,s], \calR(s,t]).
\end{align}
Since $G(x+a,y+a)=G(x,y)$ for all $x,y,a\in \bbZ^d$ we also have
\begin{align*}
G\big(\calR(r,s], \calR(s,\infty )\big)
&\stackrel{d}{=} 
G\big(\calR(0,s-r], \calR(s-r,\infty)\big)\\
&=
G\big(\calR(0,s-r]-S_{s-r},\,  \calR(s-r,\infty)-S_{s-r}\big),
\end{align*}
where $\stackrel{d}{=}$ stands for the equality in distribution. Moreover, the two random variables 
$\calR(0,s-r]-S_{s-r}$ and $\calR(s-r,\infty) - S_{s-r}$ are independent and thus we can replace the second one with an independent copy $\widetilde{\calR}(0,\infty)$. We also observe that by the symmetry of the random walk, $\calR(0,s-r]-S_{s-r}\stackrel{d}{=}\calR_{s-r-1}$.   
We infer that
\begin{align*}
\sup_{t> s}V_{r,s,t} &\leq 2 G(\calR(r,s], \calR(s,\infty)) \stackrel{d}{=} 2G(\calR(0,s-r]-S_{s-r}, \widetilde{\calR}(0,\infty))\\
&\stackrel{d}{=} 
2G(\calR_{s-r-1}, \widetilde{\calR}(0,\infty))
\leq 
2
G(\calR_{s-r}, \widetilde{\calR}_\infty ),
\end{align*}
where $\widetilde{\calR}_\infty$ is an independent copy of $\calR_\infty$.
We can easily adapt the proof of \cite[Lemma 3.2]{CSS19} to obtain that for any $p\geq 1$, 
\begin{align}\label{est-moments}
\bbE [(G(\calR_{s-r}, \widetilde{\calR}_\infty))^p]\leq c  h(s-r)^p
\end{align}
and the first inequality follows. For the second statement we proceed similarly and  \eqref{Assel-Cap-ineq} yields
\begin{align*}
\sup_{0\leq r<s}V_{r,s,t} &\leq 2 G(\calR(0,s], \calR(s,t]) 
=
2G (\calR(0,s]-S_s , \calR(s,t]-S_s)\\
&\stackrel{d}{=} 
2G (\calR_{s-1}, \widetilde{\calR}(0,t-s])
\leq 
2
G(\calR_\infty , \widetilde{\calR}_{t-s}).
\end{align*}
We combine this with \eqref{est-moments} and the proof is finished.
\end{proof}

\begin{lemma}\label{lemma_difference}
Assume \textbf{A}\ref{A1}, \textbf{A}\ref{A2} and $d/\alpha >5 /2$. For any $m\geq 1$  there is a constant $c=c(m)>0$ such that for all $s>r> 0$ \footnote{For a random variable $Y$ with finite expectation, we write $\overline{Y}=Y-\mathbb{E}[Y]$. In particular, we use the notation $\overline{\calC}_n = \calC_n - \bbE[\calC_n]$ and $\overline{\calC}(r,s] = \calC(r,s] - \bbE [\calC(r,s]]$.} 
\begin{align}\label{C-diff-estimate}
\big\Vert \overline{\calC }(0,s] - \overline{\calC} (0,r]\big\Vert_m \leq c  \sqrt{s-r}.
\end{align}
\end{lemma}

\begin{proof}
We note that by the subadditivity property of the capacity (see e.g.\ \cite[Proposition 2.2.1]{Lawler}), we have $V_{0,r,s}\geq 0$. By the second estimate of Lemma \ref{Lemma-supV} combined with Jensen's inequality we obtain
\begin{align*}
0\leq \bbE [V_{0,r,s}]^m \leq \bbE [V_{0,r,s}^m]\leq \bbE \Big[ \big( \sup_{0\leq w<r}V_{w,r,s}\big)^m \Big]  \leq c\, h(s-r)^m = o((s-r)^{m/2}),
\end{align*}
where the last asymptotics follows by the assumption $d/\alpha >5/2$. Recall that $V_{0,r,s}=\calC (0,r]+\calC(r,s]-\calC (0,s].$
This implies that we can replace $\barcalC (0,s] - \barcalC (0,r]$ with $\barcalC (r,s]$ and we can clearly assume that $r=0$. 
We next follow the approach of \cite[Lemma 3.2]{Bass_Kumagai} which is based on \cite[Lemma 4.1]{Jain-Pruitt-74}.
We observe that it suffices to show that for all $l\geq 1$,
\begin{align}\label{L-to-show}
\sup_{n\geq 1} L_{n,2l}<\infty,\quad \mathrm{where}\ \ L_{n,l}=\frac{\big\Vert \barcalC (0,n] \big\Vert_{l}}{\sqrt{n}} .
\end{align}
Indeed, we observe that by H\"{o}lder's inequality $\big\Vert \barcalC (0,n]\big\Vert_{2l-1}\leq \big\Vert \barcalC (0,n]\big\Vert_{2l}$ and thus the estimate for odd values will follow if we prove \eqref{L-to-show}. To show \eqref{L-to-show} we use induction over $l$. For $l=1$ the claim follows by the variance asymptotics from \cite[Lemmas 4.3 and 5.3]{CSS19}, see also Remark \ref{RM:Var}. We assume validity of \eqref{L-to-show} for all $l\leq l_0-1$ and we show that it holds for $l_0$. First, for any $0< r < 2n$ we have $V_{0,r,2n}\leq \sup_{t>r }V_{0,r,t}$ and thus Lemma \ref{Lemma-supV} implies that for all $l\geq 1$, $\Vert \overline{V}_{0,r,2n}\Vert_{2l}\leq c\, h(r)$. By  regular variation of $h(n)$ and by the fact that $d/\alpha >5/2$ we easily obtain
\begin{align}\label{V-o-small}
\sup_{0< r <2n} \Vert \overline{V}_{0,r,2n}\Vert_{2l}= o(\sqrt{n}).
\end{align}
This allows us to replace $\barcalC (0,2n]$ with $\barcalC (0,n ] + \barcalC (n,2n ]\stackrel{d}{=} \barcalC (0,n ] + \overline{\calC ^\prime}(0,n ]$, where $\calC ^\prime(0,n ]$ is an independent copy of $\calC (0,n ]$. 
We have
\begin{multline}\label{sum_induct}
\bbE \Big[ \big(\barcalC (0,n] + \overline{\calC^\prime}(0,n]\big)^{2l_0}\Big]
=
\sum_{j=0}^{2l_0}{2l_0 \choose j}\bbE[\overline{\calC}(0,n]^j] \bbE[\overline{\calC}(0,n]^{2l_0-j}]\\
\leq
2\bbE[\barcalC(0,n]^{2l_0}]
+
\sum_{j=2}^{2l_0 -2}
{2l_0 \choose j} \bbE \big[\left\vert \barcalC(0,n]\right\vert^j \big] \bbE \big[ \left\vert \barcalC(0,n]\right\vert ^{2l_0-j} \big].
\end{multline}
We observe that the addends for $j=1$ and $j=2l_0 -1$ disappear as $\bbE [\overline{\calC}_1]=0$. For the first term we have
\begin{align*}
2\bbE[\barcalC(0,n]^{2l_0}] = 2n^{l_0}L_{n,2l_0}^{2l_0}.
\end{align*}
Moreover, by the inductive assumption we have
\begin{align*}
\sum_{j=2}^{2l_0 -2}
{2l_0 \choose j} \bbE \big[\left\vert \barcalC(0,n]\right\vert^j \big] \bbE \big[ \left\vert \barcalC(0,n]\right\vert ^{2l_0-j} \big]
=
n^{l_0}
\sum_{j=2}^{2l_0 -2}{2l_0 \choose j}L_{n,j}^j L_{n,2l_0-j}^{2l_0-j} 
\leq c_1 n^{l_0},
\end{align*}
where $c_1=c_1(l_0) = \sup_{n,n^\prime} \sum_{j=2}^{2l_0 -2}{2l_0 \choose j}L_{n,j}^j L_{n^\prime,2l_0-j}^{2l_0-j}  <\infty $.
By combining this with \eqref{V-o-small} we obtain 
\begin{align*}
\Vert \barcalC (0,2n]\Vert_{2l_0}
\leq \sqrt{n} \big(2 L^{2l_0}_{n, 2l_0}+c_1\big)^{1/2l_0} +o(\sqrt{n}).
\end{align*}
We divide  both sides by $\sqrt{2n}$ and it follows that 
for a constant $c_2>0$,
\begin{align}\label{L-bound}
L_{2n,2l_0}\leq \left( 2^{-(l_0-1)}L_{n,2l_0}^{2l_0} +2^{-l_0}c_1 \right)^{1/2l_0} + c_2.
\end{align}
We first claim that $\{L_{2^n, 2l_0}\}_{n\geq 1}$ is bounded. To  show the claim, we choose $N$ large enough such that
\begin{align}\label{N-bound}
\left( 2^{-(l_0-1)}+ \frac{c_1}{2^{l_0}N^{2l_0}}\right)^{1/2l_0}+ \frac{c_2}{N}\leq 1.
\end{align}
Then, either $L_{k,2l_0}\leq N$ for all $k\in \bbN$ that are equal to a power of $2$ and the claim holds, or there is $k_0\in  \bbN$ which is equal to some power of $2$ such that $L_{k_0,2l_0}\geq N$. In the latter case we show by induction that 
\begin{equation}\label{induc_L}
L_{n,2l_0}\leq L_{k_0,2l_0},\quad \mathrm{for\ all}\ n\geq k_0\ \mathrm{which\ are\ equal\ to\ some\  power\ of\ 2},
\end{equation}
%$L_{n,2(l_0+1)}\leq L_{k_0,2(l_0+1)}$ for all $n\geq k_0$ that are equal to some power of $2$, 
and this will also evidently imply the claim. For $n=k_0$ \eqref{induc_L} holds trivially. Next we assume that \eqref{induc_L} holds for some $n_0\in \bbN$ which is a power of $2$ and we show that it is valid also for $2n_0$. By \eqref{L-bound},
\begin{align*}
\frac{L_{2n_0,2l_0}}{L_{k_0,2l_0}}
\leq 
\left( 2^{-(l_0-1)}\left( \frac{L_{n_0,2l_0}}{L_{k_0,2l_0}}\right)^{2l_0}+ \frac{c_1}{2^{l_0}L^{2l_0}_{k_0,2l_0}}\right)^{1/2l_0}+ \frac{c_2}{L_{k_0,2l_0}}\leq 1,
\end{align*}
where the last inequality follows by \eqref{N-bound} and the inductive assumption. Thus \eqref{induc_L} is proved and we conclude that $\{L_{2^n,2l_0}\}_{n\geq 1}$ is bounded.

We are left to show that $\{L_{n,2l_0}: n\in \bbN \text{ is not a power of } 2\}$ is bounded. To this end, we prove  that for any $n\in \bbN$ and any $m\in [2^{n-1}, 2^n)\cap\N$ it holds 
\begin{align}\label{claim_not_powers_2}
2^{-l_0}L_{m,2l_0}^{2l_0}\leq  \big(L_{2^n,2l_0}+o_n(1)\big)^{2l_0} +c_1.
\end{align} 
By the preceding part of the proof it will evidently imply the desired result. We observe that if $n=1$ then necessarily $m=1$ and \eqref{claim_not_powers_2} is true as $L_{1,2l_0}=0$. 
Indeed,
\begin{align*}
\calC(0,1]=\mathrm{Cap}(\{S_1\})=\mathrm{Cap}(\{S_1\}-S_1)=\mathrm{Cap}(\{0\})=G(0,0)^{-1}
\end{align*}
and whence $\Vert \barcalC_1\Vert_{2l_0}=0$.
For $n\geq 2$ and $m\in [2^{n-1}, 2^n)\cap\N$ we proceed as follows.
We have $\barcalC (0,m] +\overline{\calC^\prime}(0,2^n-m] \stackrel{d}{=}\barcalC(0,2^n] + \overline{V}_{0,m,2^n}$, where $\calC^\prime (0,2^n-m]$ is an independent copy of $\calC(0,2^n-m]$. Hence
\begin{align}\label{split_help}
\bbE \big[ \big( \barcalC (0,m] +\overline{\calC^\prime}(0,2^n-m] \big)^{2l_0} \big]
\leq 
\left( \big\Vert \barcalC(0,2^n] \big\Vert_{2l_0} + \big\Vert \overline{V}_{0,m,2^n}\big\Vert_{2l_0} \right)^{2l_0}.
\end{align}
Further, similarly as in \eqref{sum_induct},
\begin{multline*}
\bbE \big[ \big( \barcalC (0,m] +\overline{\calC^\prime}(0,2^n-m] \big)^{2l_0} \big]
\\=
\bbE \big[ \barcalC(0,m]^{2l_0} \big] 
+ \bbE \big[ \barcalC(0,2^n-m]^{2l_0} \big]
+
\sum_{k=2}^{2l_0-2}{2l_0 \choose k}\bbE \big[ \barcalC(0,m]^k \big]\bbE \big[ \barcalC(0,2^n-m]^{2l_0-k} \big].
\end{multline*}
This together with \eqref{split_help} leads to
\begin{multline*}
\bbE \big[ \barcalC(0,m]^{2l_0} \big] 
\\
\leq 
\left( \big\Vert \barcalC(0,2^n] \big\Vert_{2l_0} + \big\Vert \overline{V}_{0,m,2^n}\big\Vert_{2l_0} \right)^{2l_0}
+
\sum_{k=2}^{2l_0-2}{2l_0 \choose k}\bbE \big[ \big\vert\barcalC(0,m]\big\vert^k \big]\bbE \big[\big\vert \barcalC(0,2^n-m]\big\vert^{2l_0-k} \big].
\end{multline*}
Next, after dividing both sides by $2^{nl_0}$ and applying \eqref{V-o-small} we obtain
\begin{align*}
\frac{m^{l_0}}{2^{nl_0}}L_{m,2l_0}^{2l_0}
\leq \left( L_{2^n,2l_0}+o_n(1)\right)^{2l_0} + 
\sum_{k=2}^{2l_0-2}{2l_0 \choose k} L_{m,k}^kL_{2^n-m, 2l_0-k}^{2l_0-k}
\leq 
\left( L_{2^n,2l_0}+o_n(1)\right)^{2l_0} 
+c_1
\end{align*}
and this implies \eqref{claim_not_powers_2}. This finally finishes the inductive step (over $l$) and the proof is complete.
%suppose that there is $n_0\in \bbN$ such that for all  $ j\leq n_0-1$, $j\in\N$, and all   $m\in [2^{j-1}, 2^j)\cap\N$ we have $2^{-l_0}L_{m,2l_0}\leq \big(L_{2^j,2l_0}+o_j(1)\big)^{2l_0} +c_1$. Next, take any $m\in [2^{n_0-1}, 2^{n_0})\cap \bbN$. 
%I%f $m=2^{n_0}$ then we apply the previous part of the proof. For $m\in (2^{n_0}, 2^{n_0+1})\cap \bbN$
%we use the following estimate
%\begin{align*}
%\Vert \bar{\calC}_m\Vert_{2(l_0+1)}
%\leq
% \Vert \bar{\calC}_{2^{n_0+1}}\Vert_{2(l_0+1)}
% +
% \Vert \bar{\calC}_{2^{n_0+1}-m}\Vert_{2(l_0+1)}
% +
% \Vert \bar{V}_{0,m,2^{n_0+1}}\Vert_{2(l_0+1)}.
%\end{align*}
%This together with \eqref{V-o-small} yields
%\begin{align*}
%2^{-1/2}L_{m, l_0+1}\leq L_{2^{n_0+1}, l_0+1}+\frac{\sqrt{2^{n_0+1}-m}}{\sqrt{2^{n_0+1}}}L_{2^{n_0+1}-m,l_0+1}+c_6
%\end{align*} for some $c_6=c_6(l_0)>0$, 
%and the result follows by the inductive assumption. 
\end{proof}

\begin{lemma}\label{lemma_C-diff-max-estimate}
Assume \textbf{A}\ref{A1}, \textbf{A}\ref{A2} and $d/\alpha >5 /2$.
For any $m\geq 3$ there is a constant $c_2=c_2(m)>0 $ such that for all $\lambda >0$ and all $s>r> 0$
\begin{align}\label{C-diff-max-estimate}
 \bbP \big(\max_{r\leq n \leq s} \big\vert \barcalC (0,n] - \barcalC (0,r] \big\vert >\lambda \sqrt{s-r}\big)
 \leq c_2 \lambda^{-m}.
\end{align}
\end{lemma}

\begin{proof}
We adapt the approach of \cite[Lemma 3.3(b)]{Bass_Kumagai}. We set
\begin{align*}
\Theta_k &= \frac{\barcalC(0,k+r]-\barcalC(0,r]}{\sqrt{s-r}},\qquad k=0,1,2\ldots .
\end{align*}
It suffices to show that for any $m\geq 3$ there is a constant $c=c(m)>0$ such that
\begin{align*}
\bbP \big( \max_{0\le k\leq s-r}|\Theta_k|>\lambda \big)\leq \frac{c}{\lambda^m},\qquad \lambda >0.
\end{align*}
We have
\begin{align*}
\bbP \big( \max_{0\le k\leq s-r}|\Theta_k|>\lambda \big)\leq
\bbP \big( \max_{0\le k< s-r}|\Theta_k|>\lambda \big) + \bbP (|\Theta_{s-r}|>\lambda)
\end{align*}
and in view of Markov's inequality combined with \eqref{C-diff-estimate} we are left to find the desired estimate for the part with $0\le k<s-r$.
For any $0\le k <s-r$ let $k_j$ be the maximal element of \footnote{We write $\langle x \rangle = n$ if $x\in (n-\frac{1}{2}, n+\frac{1}{2}]$.} $\{ \langle \frac{i(s-r)}{2^j}\rangle : \, i=0,1,\ldots , 2^j\}$ that is smaller or equal to $k$. Evidently, there exists $J\geq 1$ such that $k_j=k$ for $j\geq J$. Thus,
\begin{align*}
\Theta_k = \Theta_{k_0} + (\Theta_{k_1}-\Theta_{k_0})+ (\Theta_{k_2}-\Theta_{k_1})+\ldots
\end{align*} 
is a finite sum. Note that 
$$
k_0 = \begin{cases}
s-r, &\quad k=s-r ,\\
0, &\quad \text{otherwise.}
\end{cases}
$$   
Hence, $|\Theta_k|> \lambda$ for some $k<s-r$ only if there exists $j\geq 1$ such that $|\Theta_{k_j}-\Theta_{k_{j-1}}|>\frac{\lambda}{2j^2}$. This in turn implies that there must exist  $i\in \{0,1,\ldots ,2^j-2\}$ such that  
\begin{align*}
k_j=\langle ( i+1)(s-r)/2^j\rangle
\quad \mathrm{and}\quad
\left\vert\Theta_{\left\langle \frac{(i+1)(s-r)}{2^j}\right\rangle} - \Theta_{\left\langle \frac{i(s-r)}{2^j}\right\rangle}\right\vert
>\frac{\lambda}{2j^2}.
\end{align*}
Let $I$ be the set of all pairs $j\in \bbN$ and $i \in \{ 0,1,\ldots , 2^j -2\}$ such that there is $k\in \{0,\ldots , s-r-1\}$ for which $k_j =\langle ( i+1)(s-r)/2^j\rangle $. 
Together with \eqref{C-diff-estimate} we infer that
\begin{align*}
\bbP \big( \max_{0\le k< s-r}|\Theta_k|>\lambda \big)
&\leq 
\sum_{(j,i)\in I}
\bbP \left( \left\vert\Theta_{\left\langle \frac{(i+1)(s-r)}{2^j}\right\rangle} - \Theta_{\left\langle \frac{i(s-r)}{2^j}\right\rangle}\right\vert
>\frac{\lambda}{2j^2} \right)\\
&\leq
\frac{c\, 2^m}{\lambda ^m(s-r)^{m/2}}\sum_{(j,i)\in I}j^{2m}\left( \left\langle \frac{(i+1)(s-r)}{2^j}\right\rangle
- \left\langle \frac{i(s-r)}{2^j}\right\rangle
\right)^{m/2}.
\end{align*}
%\leq 
%\frac{2^m}{\lambda^m}
%\sum_{j=1}^\infty j^{2m}\bbE \big[ |G_{k_j}-G_{k_{j-1}}|^m \big]\\
%&\leq 
%\frac{c_7\,  2^m}{(s-r)^{m/2}\lambda^m}
%\sum_{j=1}^{J} j^{2m}|k_j - k_{j-1}|^{m/2},
We observe that if $ (s-r)/2^j \leq 1/2$ then $\langle ( i+1)(s-r)/2^j\rangle - \langle i(s-r)/2^j\rangle =k - k =0$. Further, if $1/2 < (s-r)/2^j\leq 1$ then $\langle ( i+1)(s-r)/2^j\rangle - \langle i(s-r)/2^j\rangle =k_j-k_{j-1}\leq 1 < 2(s-r)/2^j$.
Similarly, in the case when $(s-r)/2^j>1$ we use the following estimate
\begin{align*}
k-k_{j} &\leq \left\vert k_j - \left\langle \frac{(i+2)(s-r)}{2^j}\right\rangle\right\vert -1 \\
&=
\left\langle \frac{(i+2)(s-r)}{2^j}\right\rangle
- \left\langle \frac{(i+1)(s-r)}{2^j}\right\rangle
-1
\leq \frac{s-r}{2^j}
\end{align*}
and whence $k_j - k_{j-1} \leq  k- k_{j-1}\leq 2(s-r)/2^j$. 
We thus arrive at
\begin{align*}
\bbP \big( \max_{0\le k< s-r}|\Theta_k|>\lambda \big)
&\leq 
c 2^{3m/2} \lambda^{-m}\sum_{(j,i)\in I} j^{2m}2^{-mj/2}\\
&\leq 
c 2^{3m/2} \lambda^{-m}\sum_{j=1}^{\infty}\sum_{i=0}^{2^j-2}j^{2m}2^{-mj/2}\\
&\leq 
c 2^{3m/2} \lambda^{-m}\sum_{j=1}^\infty
j^{2m}2^{(1-m/2)j}
\end{align*}
and as the last series converges (since $m>2$) the proof is finished.
\end{proof}

We proceed to the proof of Theorem \ref{TM:1} and we start with some preparation. 
We observe that 
%\begin{align*}
%\bbE [\calC(0,n_k]] = \sum_{j=1}^k \bbE [\calC(n_{j-1},n_j]] - \sum_{j=1}^{k-1}\bbE [V_{0,n_j,n_{j+1}}]
%\end{align*} 
%and thus 
centering at \eqref{Cap-decomp-Dembo-Okada} yields
\begin{align}\label{Cap-decomp-Dembo-Okada-centred}
\barcalC_{n_k}=\sum_{j=1}^k \overline{U}_j - \overline{\Delta}_{n_k,k}.
\end{align} 
Recall that the random variables $\overline{U}_j$, $j=1,\dots,k$, are independent.
%We 
%fix a parameter $\zeta >0$ which will be specified later and  
We
consider an increasing sequence $\{n_i\}_{i\ge1}\subseteq\N$ (with $n_0=0$) of all consecutive integers from each interval $[2^k,2^{k+1})$, $k\in \bbN$, which are of the form $2^k + \langle i2^k /k\rangle$, for $i=0,1,\ldots ,k-1$. We easily verify that for $n_i\in [2^k,2^{k+1})$ it holds $ n_{i+1}-n_i\leq 1+2^k/k$ and whence
\begin{align}\label{sequence_n_i}
\lim_{i\to \infty}\frac{n_{i+1}}{n_i}=1,\qquad \mathrm{and}\qquad n_{i+1}-n_i = O\left( n_i/ \log n_i \right).
\end{align}

We set  $\Lambda=\frac{d}{\alpha} - \frac{5}{2}$. According to \eqref{h_d_function}, if $\Lambda \geq 1/2$ then $h(n)$ is non-decreasing and slowly varying, see \cite[Lemma 2.2]{LeGall-Rosen}.
If $\Lambda \in (0,1/2)$ then $h(n)$ is non-decreasing and regularly varying of index $1/2-\Lambda$, which is positive and strictly smaller than $1/2$.

\begin{proof}[Proof of Theorem \ref{TM:1}]
We first show the result for the process $\{\calC (0,n]\}_{n\geq 0}$ and only in the end of the proof we switch to the original process $\{\mathrm{Cap}(\calR_n)\}_{n\geq 0}$.

Recall that $\Lambda =\frac{d}{\alpha}-\frac{5}{2}$. We fix the following constants
\begin{align}\label{constants}
\beta = 
\begin{cases}
\frac{1-\frac{3}{4}\Lambda}{2},& \Lambda \in (0,1/2);\\
\frac{5}{16},& \Lambda \geq 1/2,
\end{cases}
\quad
\varepsilon = 
\begin{cases}
\frac{\Lambda}{4},& \Lambda \in (0,1/2);\\
\frac{1}{8},& \Lambda \geq 1/2,
\end{cases}
\quad
\gamma = \frac{1}{4}.
\end{align}
We start by handling the error term in \eqref{Cap-decomp-Dembo-Okada-centred} and prove that it is negligible, that is, there exists a constant $c_1>0$ such that 
\begin{align}\label{error-to-zero}
\limsup_{i\to \infty}\frac{\left\vert \overline{\Delta}_{n_i,i}\right\vert}{n_i^\beta}\leq c_1.
\end{align}
We set $V_j  = V_{0,n_j,n_{j+1}}$. By H\"{o}lder's inequality combined with the second inequality in Lemma \ref{Lemma-supV}, 
\begin{align*}
\bbE \Big[ \Big( \sum_{j=1}^{i-1}V_j\Big)^4\Big] 
&=
\sum_{j_1,\ldots ,j_4=1}^{i-1}\bbE [V_{j_1}\cdots V_{j_4}]
\leq 
\sum_{j_1,\ldots ,j_4=1}^{i-1} \left( \prod_{r=1}^4 \bbE \big[ V_{j_r}^4\big] \right)^{1/4}\\
&\leq 
c_2 \sum_{j_1,\ldots ,j_4=1}^{i-1} \prod_{r=1}^4 h(n_{j_r+1} - n_{j_r})
=
c_2 \left( \, \sum_{j=1}^{i-1}h(n_{j+1}-n_j)\right)^4 ,
\end{align*} 
for some $c_2>0$.
Thus, for any $k_0\in \bbN$ and $2^{k_0}\leq n_i <2^{k_0+1}$ we obtain 
\begin{align}\label{Delta_est_help}
\bbP \left( \Delta_{n_i,i}\geq c_1n_i^\beta \right) 
&\leq \frac{\bbE \Big[ \Big( \sum_{j=1}^{i-1}V_j\Big)^4\Big]}{c_1^4n_i^{4\beta}}
\leq c_3 2^{-4\beta k_0}\left( \sum_{k=1}^{k_0}k\,   h(2^k/k )\right)^4.
\end{align}
We next distinguish between two cases. If $\Lambda \in (0,1/2)$ then we have $h(n) \leq c_4 n^{\frac{1}{2}-\Lambda+\varepsilon}$, for some $c_4>0$ and all  $n\in\N$. 
Hence
\begin{align*}
\bbP \left( \Delta_{n_i,i}\geq c_1n_i^\beta \right) 
\leq 
c_5 2^{-4\beta k_0} \left( \sum_{k=1}^{k_0}k^{\frac{1}{2}+\Lambda -\varepsilon }\,  2^{k(\frac{1}{2}-\Lambda +\varepsilon)} \right)^4 ,
\end{align*} 
for some $c_5>0$.
Further, the following asymptotics can be derived via the Stolz-Ces\`{a}ro theorem:  for all $p\in \bbR$ and $q>1$, 
\begin{align}\label{aymp_series}
\lim_{k_0 \to \infty} k_0^{-p}q^{-k_0}\sum_{k=1}^{k_0}k^p q^k = \frac{q}{q-1} .
\end{align}
Hence, for some $c_6>0$,
\begin{align*}
\bbP \left( \Delta_{n_i,i}\geq c_1n_i^\beta \right) 
\leq 
c_6\, 
k_0^{2+4\Lambda-4\varepsilon} \, 2^{4k_0 (\frac{1}{2}-\Lambda+\varepsilon - \beta)}.
\end{align*}
Since the number of $n_i$'s in each interval $[2^{k_0},2^{k_0+1})$ is equal to $k_0$, we obtain that for a constant $c_7>0$,
\begin{align*}
\sum_{i=1}^\infty \bbP \left( \Delta_{n_i,i}\geq c_1n_i^\beta \right) 
\leq c_7
\sum_{k_0=1}^\infty  
k_0^{3+4\Lambda-4\varepsilon} \, 2^{4k_0 (\frac{1}{2}-\Lambda+\varepsilon - \beta)} <\infty ,
\end{align*}
where convergence of the last series follows by our choice of $\varepsilon$ and $\beta$ from \eqref{constants} as  $\frac{1}{2}-\Lambda +\varepsilon <\frac{1}{2}$ and $\beta \in (\frac{1}{2}-\Lambda +\varepsilon , \frac{1}{2})$.
If $\Lambda \geq 1/2$ then there is $c>0$ such that $h(n)\leq c n^\varepsilon
$, for all $n\in \bbN$. Using the same reasoning as above and combining \eqref{Delta_est_help} and \eqref{aymp_series} we obtain that for a constant $c_8>0$,
\begin{align*}
\sum_{i=1}^\infty \bbP \left( \Delta_{n_i,i}\geq c_1n_i^\beta \right) 
\leq c_8
\sum_{k_0=1}^\infty  
k_0^{5-4\varepsilon} \, 2^{4k_0 (\varepsilon - \beta)}
\end{align*}
and the last series converges by \eqref{constants} as $\varepsilon < \beta$.
Finally, by the Borel-Cantelli lemma we infer that
\begin{align*}
\limsup_{i\to \infty} \frac{\Delta_{n_i,i}}{n_i^\beta} \leq c_1.
\end{align*}
We proceed similarly for the expectations. There are constants $c_9,c_{10}, c_{11}>0$ such that
for any $k_0\in \bbN$ and $2^{k_0}\leq n_i <2^{k_0+1}$ it holds
\begin{align*}
\bbE [\Delta_{n_i,i}]
=
\sum_{j=1}^{i-1}\bbE [V_j]
&\leq c_9\sum_{k=1}^{k_0} k\, h(2^k/k)\\
&\leq c_{10}
\begin{cases}
k_0^{ \frac{1}{2}+\Lambda -\varepsilon }2^{k_0(\frac{1}{2}-\Lambda +\varepsilon)},& \Lambda \in (0,1/2);\\
k_0^{1-\varepsilon }2^{k_0\varepsilon},& \Lambda \geq 1/2
\end{cases}\\
&\leq 
c_{11}
\begin{cases}
\big( \log n_i\big)^{\frac{1}{2}+\Lambda -\varepsilon }\, n_i^{\frac{1}{2}-\Lambda+\varepsilon},& \Lambda \in (0,1/2);\\
\big( \log n_i \big)^{1-\varepsilon }n_i^{\varepsilon},& \Lambda \geq 1/2 
\end{cases}\\
&= o(n_i^\beta),
\end{align*} 
where the last equality follows by our choice of $\beta$ and $\varepsilon$ from \eqref{constants}. This implies \eqref{error-to-zero}.

Let $H_j = \sigma ^{-1}\overline{U}_j$, $j=1,\dots,k$, where $\sigma$ is the constant from \eqref{Var-asymp}. Recall that these random variables are independent. We apply the  Skorohod embedding theorem (see \cite{Skorohod}) which states that  there 
exist a standard Brownian motion $\{B_u\}_{u\ge0}$  and  non-negative independent stopping times $T_1, T_2,\dots$ such that for every  $k\geq 1$
\begin{align}\label{Skorohod-law}
\{B_{T_1+T_2+\ldots +T_j}\}_{j=1}^k\stackrel{d}{=} \{H_1+\ldots +H_j\}_{j=1}^k.
\end{align}
We can choose the probability space in such a way that $\{B_u\}_{u\ge0}$ and $\{S_n\}_{n\geq 0}$ are defined on the same $(\Omega,\calF,\bbP)$, see \cite[Section 3.5]{Phillip_Stout}.
Moreover, the following moment estimates hold
\begin{align}\label{T_j-moments}
\mathbb{E}[T_j]=\bbE [H_j^2]\quad \mathrm{and}\quad
\mathbb{E}[T_j^m]\le c_{12}\, \mathbb{E}[H_j^{2m}],\qquad m\geq 2,
\end{align}
for some $c_{12}=c_{12}(m)>0$. To get rid of the randomness coming from the sequence of stopping times $T_j$ in \eqref{Skorohod-law}, we first show that 
\begin{align}\label{T_j-centred-limsup}
\limsup_{i\to \infty}\frac{\sum_{j=1}^i \overline{T}_j}{n_i (\log n_i)^{-\gamma}}<\infty\qquad \mathrm{a.s.}
\end{align}
As $\sum_{j=1}^i \overline{T}_j$ is a martingale, we obtain by Doob's submartingale inequality that for any $l\geq 1$,
\begin{align}\label{T_j-Doob}
\bbP \Big( \sup_{r\leq i}\Big\vert \sum_{j=1}^r \overline{T}_j\Big\vert \geq n_i (\log n_i)^{-\gamma}\Big)
\leq
\frac{\bbE \big[ \big(\sum_{j=1}^i \overline{T}_j \big) ^{2l}\big]}{n_i^{2l}(\log n_i)^{-2\gamma l}}.
\end{align}
To estimate the right-hand side of \eqref{T_j-Doob} we find that
\begin{align}\label{T_j-combinatorial}
\begin{split}
\bbE \big[ \big( \sum_{j=1}^i \overline{T}_j \big) ^{2l}\big]
&=
\sum_{j_1,\ldots ,j_{2l}=1}^i \bbE [\overline{T}_{j_1}\cdots \overline{T}_{j_{2l}}]\\
&=
\sum_{I}\frac{(2l)!}{a_1!\cdots a_p!}\sum_{j_1,\ldots ,j_{p}=1}^i  \bbE \big[ \overline{T}_{j_1}^{a_1}\big]\cdots \bbE \big[ \overline{T}_{j_p}^{a_p}\big],
\end{split}
\end{align}
where $I$ ranges over all $(a_1,\ldots ,a_p)$ and $1\leq p\leq 2l$ such that $a_1,\ldots ,a_p \geq 2$ and $\sum_{\nu =1}^p a_{\nu}=2l$ (which evidently implies that $p\leq l$). The second equality in \eqref{T_j-combinatorial} follows from the fact that $\bbE [\overline{T}_{j_1}\cdots \overline{T}_{j_{2l}}]=0$ if there is one $j_k$ which is different from all the others as $\overline{T}_{j}$ are independent with $\bbE [\overline{T}_{j}]=0$.
Further, by \eqref{T_j-moments} and \eqref{C-diff-estimate}, for any $m\geq 1$ there exist  constant $c_{13}=c_{13}(m)>0$ and $c_{14}=c_{14}(m)>0$ such that
\begin{align*}
\big\vert \bbE \big[ \overline{T}_j^m \big] \big\vert
&=
\Big\vert \bbE \Big[ \sum_{k=0}^{m}{m\choose k}T_j^k (-1)^{m-k}\big( \bbE [T_j]\big)^{m-k} \Big]  \Big\vert \\
&\leq
c_{13}
\sum_{k=0}^{m}{m\choose k}\bbE \big[ H_j^{2k} \big]\big( \bbE [H_j^2] \big)^{m-k}
\leq c_{14}(n_{j} - n_{j-1})^m.
\end{align*}
Hence, by \eqref{sequence_n_i} and  \eqref{aymp_series} we have for $2^{k_0}\leq n_i < 2^{k_0+1}$,
\begin{align*}
\sum_{j=1}^i \big\vert \bbE \big[ \overline{T}_j^m\big] \big\vert 
\leq
 c_{14}\sum_{j=1}^i (n_{j}- n_{j-1})^m\leq c_{15}\sum_{k=1}^{k_0}k^{1-m}2^{mk}
  \leq
 c_{16}\, 
 k_0^{1-m}2^{mk_0} ,
\end{align*} 
for some $c_{15},c_{16}>0.$
Applying this to \eqref{T_j-combinatorial} we conclude that for any $l\geq 1$ there exist some constants $c_{17}=c_{17}(l)>0$ and $c_{18}=c_{18}(l)>0$ such that
\begin{equation}
\begin{aligned}\label{T_j-2l-moment}
\bbE \Big[ \Big\vert \sum_{j=1}^i \overline{T}_j \Big\vert ^{2l}\Big]
&\leq
c_{17}\sum_{I}
\frac{(2l)!}{a_1!\cdots a_p!}
\prod_{\nu=1}^p k_0^{1-a_\nu}2^{a_\nu k_0}
\\&=
c_{17}k_0^{-2l}2^{2lk_0}
\sum_{I}
\frac{(2l)!}{a_1!\cdots a_p!} k_0^{ p}
\leq
 c_{18} k_0^{-l} 2^{2lk_0},
\end{aligned}
\end{equation}
where the last inequality follows as $p\leq l$. Coming back to \eqref{T_j-Doob} yields
\begin{align*}
\sum_{i=1}^{\infty} \bbP \Big( \sup_{r\leq i}\Big\vert \sum_{j=1}^r \overline{T}_j\Big\vert \geq n_i (\log n_i)^{-\gamma}\Big)
\leq 
c_{19}
\sum_{k_0=1}^\infty k_0 \frac{k_0^{-l}2^{2lk_0}}{2^{2lk_0}k_0^{-2\gamma l}}
=c_{19}
\sum_{k_0=1}^\infty k_0^{1-l(1-2\gamma)}.
\end{align*}
By \eqref{constants} we know that $\gamma <1/2$ and thus we can choose $l$ large enough such that $1-l(1-2\gamma) <-1$. Then the last series converges and \eqref{T_j-centred-limsup} follows by the Borel-Cantelli lemma.

Further, let $J_i=\sum_{j=1}^iT_j$, $\delta_i = n_i (\log n_i)^{-\gamma}$. We clearly have
\begin{multline}\label{B-est-1}
\sum_{i=1}^\infty \bbP \big( |B_{J_i}-B_{\bbE[J_i]}|>\delta_i^{1/2}(\log n_i)^{\gamma /2} \big) \\
\leq 
 \sum_{i=1}^\infty \bbP \big( |B_{J_i}-B_{\bbE[J_i]}|>\delta_i^{1/2}(\log n_i)^{\gamma /2},\, |\overline{J}_i|\leq 2\delta_i \big)
 +
 \sum_{i=1}^\infty \bbP \big( |\overline{J}_i|>2\delta_i \big).
\end{multline}
We first handle the first term of the right-hand-side of \eqref{B-est-1}. Let
\begin{equation*}
	\hat{B}(r,s)= \sup_{r\le u, v\le s}|B_v-B_u|.
	\end{equation*} 
	We claim that
	\begin{equation*}
	\hat{B}(r,s)\stackrel{\text{(d)}}{=} \sup_{0\le u, v\le s-r}|B_v-B_u|\le 2\sup_{0\le u\le s-r}|B_u|.
	\end{equation*}
	Indeed, for any $w\ge0$ we have 
	\begin{align*}
	\Prob(\hat{B}(r,s)\le w)&=\mathbb{E}[\Ind_{\{\hat{B}(r,s)\le w\}}]= \mathbb{E}[\mathbb{E}[\Ind_{\{\hat{B}(0,s-r)\le w\}}\circ\theta_r|\mathcal{F}_r]]\\&=\int_{\R}\mathbb{E}[\Ind_{\{\hat{B}(0,s-r)\le w\}}]\Prob_{B_r}(\D x)=\Prob(\hat{B}(0,s-r)\le w),
	\end{align*} where in the third step we used Markov property and space homogeneity of $\{B_u\}_{u\ge0}$, $\{\theta_u\}_{u\ge0}$ denotes the standard shift operator, and $\Prob_{B_r}(\D x)$ stands for the distribution of the random variable $B_r$.
This implies
\begin{multline}\label{sup-B-est}
\bbP \big( |B_{J_i}-B_{\bbE[J_i]}|>\delta_i^{1/2}(\log n_i)^{\gamma /2},\, |\overline{J}_i|\leq 2\delta_i \big)\\
\leq
\bbP \Big( \sup_{\bbE[J_i]-2\delta_i\leq u,v\leq \bbE[J_i]+2\delta_i}|B_v-B_u|> \delta_i^{1/2}(\log n_i)^{\gamma /2} \Big)\\
\leq
\bbP \Big( \sup_{0\leq u\leq 4\delta_i}|B_u|>\frac{1}{2} \delta_i^{1/2}(\log n_i)^{\gamma /2} \Big)
\leq
2e^{  -\frac{(\log n_i)^{\gamma }}{32}},
\end{multline}
	where in the last step we applied \cite[Excercise II.1.23]{RY_book}.
We obtain that for $2^{k_0}\leq n_i <2^{k_0+1}$ and for some constants $c_{20}, c_{21}>0$,
\begin{align*}
 \sum_{i=1}^\infty \bbP \big( |B_{J_i}-B_{\bbE[J_i]}|>\delta_i^{1/2}(\log n_i)^{\gamma /2},\, |\overline{J}_i|\leq 2\delta_i \big)
 \leq 
 c_{20}
 \sum_{k_0=1}^\infty 
 k_0 e^{-c_{21} k_0^{\gamma }}<\infty.
\end{align*}
For the second term in \eqref{B-est-1} we apply \eqref{T_j-2l-moment} and obtain
\begin{align*}
\sum_{i=1}^\infty \bbP \big( |\overline{J}_i|>2\delta_i \big)
&\leq 
\sum_{i=1}^{\infty}
\frac{\bbE \Big[ \big\vert \sum_{j=1}^i \overline{T}_j \big\vert^{2l} \Big]}{(2\delta_i)^{2l}}\\
&\leq
c_{18}\sum_{k_0=1}^\infty k_0
\frac{k_0^{-l}2^{2lk_0}}{2^{2l}2^{2lk_0}k_0^{-2\gamma l}}
=
c_{18}2^{-2l}
\sum_{k_0=1}^\infty k_0^{1-l(1-2\gamma)}<\infty .                                                                                                                                                                                                                                                                                                                                                                                                   
\end{align*} 
Together with \eqref{B-est-1} we arrive at
\begin{align}\label{B-est-2}
B_{J_i} - B_{\bbE[J_i]} = O(\delta_i^{1/2} (\log n_i)^{\gamma /2}) = O(\sqrt{n_i})\qquad \mathrm{a.s.}
\end{align}
Next, by Lemma \ref{LM:3.2} and Remark \ref{RM:Var} there is a constant $c_{22}>0$ such that
\begin{align*}
\bbE[J_i] 
&=
\sum_{j=1}^i \bbE [H_j^2] = \sigma^{-2}\sum_{j=1}^i \Var \big(\calC (0,n_j-n_{j-1}]\big)\\
& \leq \sum_{j=1}^i \big( (n_j-n_{j-1})+ c_{22} \sqrt{n_j-n_{j-1}}h(n_j-n_{j-1})\big)
\end{align*}
and similarly as before we obtain for some constants $c_{23},c_{24}, c_{25}, c_{26}>0$,
\begin{align*}
\sum_{j=1}^i  \sqrt{n_j-n_{j-1}}h(n_j-n_{j-1})
&\leq
c_{23}
\sum_{k=1}^{k_0} k \sqrt{2^k/k}\, h(2^k/k)\\
&\leq
c_{24} 
\begin{cases}
\sum_{k=1}^{k_0}k^{\Lambda -\varepsilon}2^{k(1-\Lambda +\varepsilon)},& \Lambda \in (0,1/2);\\
\sum_{k=1}^{k_0} k^{\frac{1}{2}-\varepsilon}2^{k(\frac{1}{2}+\varepsilon)},& \Lambda \geq 1/2
\end{cases}\\
&\leq c_{25}
\begin{cases}
k_0^{\Lambda -\varepsilon}2^{k_0 (1-\Lambda +\varepsilon)},& \Lambda \in (0,1/2);\\
k_0^{\frac{1}{2}-\varepsilon}2^{k_0(\frac{1}{2}+\varepsilon)},& \Lambda \geq 1/2
\end{cases}\\
&\leq c_{26}
\begin{cases}
n_i^{1-\Lambda +\varepsilon}(\log n_i)^{\Lambda-\varepsilon},& \Lambda \in (0,1/2);\\
n_i^{\frac{1}{2}+\varepsilon}(\log n_i)^{\frac{1}{2}-\varepsilon},& \Lambda \geq 1/2
\end{cases}\\
&=
\begin{cases}
o(n_i^{1-\Lambda + 2\varepsilon}),& \Lambda \in (0,1/2);\\
o(n_i^{\frac{1}{2}+2\varepsilon}),& \Lambda \geq 1/2,
\end{cases}
\end{align*} 
where  the last equality is justified by  the fact that 
the log-terms can be killed by $n_i^{\varepsilon /2}$ in view of the slow variation. We set $v_i = \bbE[J_i]$
and we infer that
\begin{align}\label{J-est}
v_i \leq 
\begin{cases}
n_i + o\big(n_i^{1-\Lambda+2\varepsilon}\big),& \Lambda \in (0,1/2);\\
n_i + o\big( n_i^{\frac{1}{2}+2\varepsilon}\big),& \Lambda \geq 1/2.
\end{cases}
\end{align}
Note that by \eqref{constants} the exponents in the error terms in \eqref{J-est} are both smaller than one.
We next claim that (in both cases: $\Lambda \in (0,1/2)$ and $\Lambda \geq 1/2$)
\begin{align}\label{B-est-3}
\left\vert B_{v_i}-B_{n_i}  \right\vert=O(\sqrt{n_i})\qquad \mathrm{a.s.}
\end{align}
Indeed, similarly as in \eqref{sup-B-est}, there exist positive constants $ c_{27}, c_{28}, c_{29}, c_{30}, c_{31}$ such that for $2^{k_0}\leq n_i < 2^{k_0+1}$, 
\begin{align*}
\sum_{i=1}^\infty
\bbP \Big( \left\vert B_{v_i}-B_{n_i}  \right\vert > \sqrt{n_i} \Big)
&\leq 2
\begin{cases}
\sum_{i=1}^\infty e^{ -c_{27}n_i^{\Lambda -2\varepsilon}},& \Lambda \in (0,1/2);\\
\sum_{i=1}^\infty e^{-c_{28}n_i^{\frac{1}{2}-2\varepsilon}},& \Lambda \geq 1/2
\end{cases}\\
&\leq c_{29}
\begin{cases}
\sum_{k_0=1}^\infty k_0\, e^{ -c_{30}2^{k_0(\Lambda -2\varepsilon)}},& \Lambda \in (0,1/2);\\
\sum_{k_0=1}^\infty k_0\,  e^{-c_{31}2^{k_0(\frac{1}{2}-2\varepsilon)}},& \Lambda \geq 1/2
\end{cases}
<\infty,
\end{align*}
where convergence of the series is justified by our choice of $\varepsilon$ from« \eqref{constants}.
Combining \eqref{B-est-2} with \eqref{J-est} and \eqref{B-est-3} we arrive at
\begin{align*}
B_{T_1+\ldots + T_i}-B_{n_i}=O(\sqrt{n_i})\qquad \mathrm{a.s.}
\end{align*}
Taking into account all the obtained asymptotics we conclude that
\begin{align*}
\sigma^{-1}\barcalC(0,n_i]-B_{n_i}
&=
\sigma^{-1} \sum_{j=1}^i \overline{U}_j + \sigma^{-1}\overline{\Delta}_{n_i,i}-B_{n_i}\\
&= \sigma^{-1} \sum_{j=1}^i \overline{U}_j -B_{n_i} + O(n_i^\beta)\quad 
\mathrm{a.s.}
\end{align*}
Further, as $\sigma^{-1} \sum_{j=1}^i \overline{U}_j $ has the same law as $B_{T_1+\ldots +T_i}$ (see \eqref{Skorohod-law}) and they are defined on the same probability space, we are allowed to replace the first term in the asymptotics above and we obtain
\begin{align}\label{final_asymp}
\sigma^{-1}\barcalC(0,n_i]-B_{n_i}
&=
B_{T_1+\ldots +T_i§}-B_{n_i}+O(n_i^\beta) = O(\sqrt{n_i})\quad \mathrm{a.s.}
\end{align}

We next show that \eqref{final_asymp} is valid for all $n\in\N$. By \eqref{C-diff-max-estimate}, for $2^{k_0}\leq n_i <2^{k_0+1}$,
\begin{align*}
\sum_{i=1}^\infty \bbP \Big( \max_{n_i\leq n\leq n_{i+1}}|\barcalC(0,n] - \barcalC(0,n_i]| > \sqrt{n_i}\Big) 
\leq 
c_{32}\sum_{i=1}^\infty \Big( \frac{n_{i+1}-n_i}{n_i}\Big)^{m/2}
\leq
c_{33}\sum_{k_0=1}^\infty k_0^{1-m/2},
\end{align*} 
for some $c_{32},c_{33}>0$,
and the last series converges 
if we choose $m >4$. 
Hence 
\begin{align*}
\max_{n_i\leq n\leq n_{i+1}}|\barcalC(0,n] - \barcalC(0,n_i]| = O(\sqrt{n_i})\qquad \mathrm{a.s.}
\end{align*}
Similarly,  for $2^{k_0}\leq n_i <2^{k_0+1}$,
\begin{align*}
\sum_{i=1}^\infty \bbP \Big( \sup_{n_i\leq u\leq n_{i+1}}|B_u -B_{n_i}|>\sqrt{n_i}\Big) 
&\leq
2
\sum_{i=1}^\infty
e^{- n_i/8(n_{i+1}-n_i)}
\leq
c_{34}
\sum_{k_0=1}^\infty k_0\, e^{-c_{35}\, k_0},
\end{align*}
for some $c_{34}, c_{35}>0$, 
and the last series converges as well. We conclude that 
\begin{align*}
\sup_{n_i\leq u\leq n_{i+1}}|B_u -B_{n_i}| = O(\sqrt{n_i})\quad \mathrm{a.s.}
\end{align*} 
and we infer that \eqref{final_asymp} holds for all $n\in \bbN$.

As it was mentioned at the beginning of the proof, we finally present how to switch from $\calC(0,n]$ to $\calC_n=\mathrm{Cap}(\calR_n)$. 
According to \cite[Proposition 2.2.1]{Lawler} we have 
\begin{align*}
0\le 
\calC_n - \calC (0,n]
\le {\rm Cap}(\{S_0\})=1/G(0,0),
\end{align*} 
which clearly implies that
$$\big\vert \barcalC_n-  \barcalC(0,n]\big\vert
\le 
\calC_n-\calC(0,n]
+
\bbE[\calC_n]-\bbE[\calC(0,n]]
\le 
2/G(0,0).
$$
Hence, we can replace $\barcalC (0,n]$ with $\barcalC_n$ in \eqref{final_asymp} (for arbitrary $n$)  and the proof is complete. 
\end{proof}

\section{The Cardinality Process}\label{Sec:cardinality}

The goal of this section is to prove  \Cref{TM:2}. 
 Results in this section correspond to the results from  \Cref{Sec:capacity} and the main arguments can be easily repeated, thus we only formulate the corresponding results and present the main steps of the proofs. 
 For any $r,s,t \in \bbZ$, $0\leq r< s< t$ we use the notation 
 \begin{align*}
  \calR(r,s] = \{S_{r+1}, \ldots , S_{s}\}\quad\text{and}\quad V_{r,s,t}=\calR(r,s]\cap \calR(s,t].
 \end{align*}
 By translation invariance $\calR(r,s]$ has the same law as $\calR(0,s-r]$,  by the independence of the increments    $\calR(r,s]$ and $\calR(s,t]$ are independent, and $0\le|\calR_s|-|\calR(0,s]|\le1$. For any  increasing sequence $\{n_k\}_{k\ge0}\subset \bbZ$ with $n_0=0$ we have 
 \begin{align*}
 |\calR(0,n_k] |
 =
 \sum_{j=1}^k U_j - \Delta_{n_k,k},
 \end{align*}
 where $U_j = |\calR(n_{j-1},n_j]|$, $j=1,\dots,k$, are independent   and
 \begin{align*}
 \Delta_{n_k,k} = 
 \sum_{j=1}^{k-1} |V_{0,n_j,n_{j+1}}|.
 \end{align*}
 Further, let 
\begin{equation}\label{Function_H}
 H(n)=\left\{
 \begin{array}{ll}
 1, & d/\alpha>2;\\
 \sum_{k=1}^nk^{-1}(\ell(k))^{-d}, & d/\alpha=2;\\
 n^{2-d/\alpha} (\ell (n))^{-d},& 1< d/\alpha <2,
 \end{array}
 \right.
 \end{equation}
where $\ell(n)$ is the slowly varying function  from  \eqref{b_function}. Note that $H$ is slowly varying if $d/\alpha \geq 2$, see \cite[Lemma 2.2]{LeGall-Rosen}, and regularly varying of index $2-d/\alpha$ if $1<d/\alpha <2$. We first derive a version of Lemma \ref{Lemma-supV} for the cardinality process.

 \begin{lemma}\label{Lemma2-supV}
 	Assume \textbf{A}\ref{A1} and $d/\alpha >1$.
 	 For any $p\geq 1$ there is a constant $c=c(d,\alpha ,p)>0$ such that for any integers $0\leq r < s < t$,
 	\begin{align*}
 	\Vert \sup_{t>s} |V_{r,s,t}|\Vert_{p} \leq c H(s-r)\qquad \mathrm{and}\qquad 
 	\Vert \sup_{0\leq r < s} |V_{r,s,t}|\Vert_{p} \leq c H(t-s).
 	\end{align*}
 \end{lemma}
 
 \begin{proof}
% 	By translation invariance, we can assume that $r=0$. 
 	We evidently have 
 	\begin{align*}
 	\sup_{t>s} |V_{r,s,t}| & =|\calR(r,s] \cap \calR (s,\infty)| = |(\calR(r,s]-S_s) \cap (\calR (s,\infty)-S_s)|\\
 	&\stackrel{d}{=}
 	|(\calR(r,s]-S_s) \cap \widetilde{\calR} (0,\infty)|  \stackrel{d}{=}
 	|\calR_{s-r-1}\cap \widetilde{\calR} (0,\infty)| \\
 	&\leq
 	|\calR_{s-r}\cap \widetilde{\calR}_\infty|,
 	\end{align*}
 	where $\widetilde{\calR}(0,\infty)$ is an independent copy of $\calR (0,\infty)$. By
 	\cite[Corollary 3.2 and Remark on p.\ 666]{LeGall-Rosen}, 
\begin{align*}
\||\calR_n\cap\widetilde{\calR}_n|\|_p\le cH(n),
\end{align*} 	
 	 for some $c=c(d,\alpha,p)>0$.
It is straightforward to adapt the arguments of the above statement to extend it to the case when one of the indices is set to be equal to infinity and thus the first inequality follows. The second statement can be proved analogously. 
 \end{proof}

 The proof of versions of Lemma \ref{lemma_difference} and Lemma \ref{lemma_C-diff-max-estimate} for the cardinality process are analogous to those for the capacity process and thus we only state the result.

 \begin{lemma}\label{lemma_difference2}
 	Assume \textbf{A}\ref{A1} and $d/\alpha >3 /2$. For any $m\geq 1$  there is a constant $c_1=c_1(m)>0$ such that for all $s>r\geq 0$
 	\begin{align*}
 	\big\Vert \overline{|\calR_s|} - \overline{|\calR_r|}\big\Vert_m \leq c_1 \sqrt{s-r}.
 	\end{align*}
 	Further, for any $m\geq 3$ there is a constant $c_2=c_2(m)>0 $ such that for all $\lambda >0$ and all $s>r\geq 0$
 	\begin{align*}
 	\bbP \big(\max_{r\leq n \leq s} \big\vert \overline{|\calR_n|}-\overline{|\calR_r|}\big\vert>\lambda \sqrt{s-r}\big)
 	\leq c_2\, \lambda^{-m}.
 	\end{align*}
 \end{lemma}

Finally, using  the same arguments as in the proof of Theorem \ref{TM:1}, with the only difference that in this case we use the variance asymptotics of the cardinality process given in Lemma \ref{lm:var_asymp-card}, we conclude Theorem \ref{TM:2}.

\section{Appendix: Variance Second Order Asymptotics}\label{Appendix}

\subsection*{The capacity process}
We derive a second order asymptotics for the variance of the random variable $\calC_n=\mathrm{Cap}(\calR_n)$. 
According to \cite[Proposition 2.11]{Spitzer} and \cite[Proposition 1.2]{Asselah_Zd}, for $A,B\subseteq\ZZ^d$, it holds that
\begin{align*}
{\rm Cap}(A)+{\rm Cap}(B)-2G(A,B)\le {\rm Cap}(A\cup B)\le  {\rm Cap}(A)+{\rm Cap}(B)-{\rm Cap}(A\cap B).
\end{align*} 
Since 
\begin{align*}
\calC_{n+m}&=
{\rm Cap}(\calR_{n + m}-S_n)\\
&={\rm Cap}(\{S_0-S_n,\dots,S_n-S_n\}\cup\{S_n-S_n,\dots,S_{n+m}-S_n\}),
\end{align*}
we conclude from the previous relation that
\begin{equation}\label{E:cap_decomp}
\calC_n^{(1)}+\calC_m^{(2)}-2\mathcal{E}(n,m)\le \calC_{n+m}\le  \calC_n^{(1)}+\calC_m^{(2)},
\end{equation} 
where $\calC_n^{(1)}={\rm Cap}(\calR^{(1)}_n)$ and $\calC_m^{(2)}={\rm Cap}(\calR^{(2)}_m)$ with $\calR_n^{(1)}$ and $\calR_m^{(2)}$ being independent and having the same law as $\calR_n$ and $\calR_m$ respectively, and  
$	\mathcal{E}(n,m)$ has the same law as $G(\calR^{(1)}_{n}, \calR^{(2)}_{m})$.

Recall the definition of the function $h(n)$ from \eqref{h_d_function} and that $\Lambda=\frac{d}{\alpha} - \frac{5}{2}$.
%We fix the constant $\Lambda=\frac{d}{\alpha} - \frac{5}{2}$, and observe that if $\Lambda \geq 1/2$ then $h(n)$ is non-decreasing and slowly varying.
If $\Lambda \in (0,1/2)$ we present $h(n)$ in the form 
\begin{equation}\label{s_calC}
h(n)\,=\,n^{1/2-\Lambda}s(n)\,,
\end{equation} 
with $s(n)=(\ell(n))^{-d}$ being slowly varying. 	
%Clearly, in this case, $h(n)$ is then non-decreasing and regularly varying of index $1/2-\Lambda$, which is positive and  strictly smaller than $1/2$.

\begin{lemma}\label{LM:3.2} 
	Assume \textbf{A}\ref{A1}, \textbf{A}\ref{A2} and $d/\alpha >5/2$.
	It holds  	\begin{equation*}
\Var(\calC_n)= 
	\sigma^2 n+O(\sqrt{n}\, h(n)).
	\end{equation*} 
\end{lemma}
\begin{proof}
	By taking expectation in \eqref{E:cap_decomp} and then by subtracting those two relations, we get
	\begin{equation*} 
	\overline{\calC}_n^{(1)}+\overline{\calC}_m^{(2)}-2\mathcal{E}(n,m)\le \overline{\calC}_{n+m}\le  \overline{\calC}_n^{(1)}+\overline{\calC}_m^{(2)}+2\mathbb{E}[\mathcal{E}(n,m)],\end{equation*} 
	which implies
	\begin{align*}
	|\overline{\calC}_{n+m}-(\overline{\calC}_n^{(1)}+\overline{\calC}_m^{(2)})|
	\le 2\bigl(\mathcal{E}(n,m)+\mathbb{E}[\mathcal{E}(n,m)]\bigr)
	\le 2\bigl(\mathcal{E}(n,m)+\lVert \mathcal{E}(n,m)\rVert_2\bigr).
	\end{align*}
	Thus
	\begin{equation*}
	\lVert \overline{\calC}_{n+m}-(\overline{\calC}_n^{(1)}+\overline{\calC}_m^{(2)})\rVert_2\le 4\lVert \mathcal{E}(n,m)\rVert_2.
	\end{equation*}
	From this we obtain
	\begin{align*}
	\lVert\overline\calC_{n+m}\rVert_2&\le \lVert \overline{\calC}_n^{(1)}+\overline{\calC}_m^{(2)}\rVert_2+ \lVert \overline{\calC}_{n+m}-(\overline{\calC}_n^{(1)}+\overline{\calC}_m^{(2)})\rVert_2\\
	&\le \lVert \overline{\calC}_n^{(1)}+\overline{\calC}_m^{(2)}\rVert_2 + 4\lVert \mathcal{E}(n,m)\rVert_2\\
	&= \bigl(\lVert \overline{\calC}_n^{(1)}+\overline{\calC}_m^{(2)}\rVert^2_2 \bigr)^{1/2}+ 4\lVert \mathcal{E}(n,m)\rVert_2.
	\end{align*}
	By independence of $\calC_n^{(1)}$ and $\calC_m^{(2)}$, we conclude that
	\begin{equation*}
	\lVert\overline\calC_{n+m}\rVert_2\le\bigl(\lVert \overline{\calC}_n^{(1)}\rVert_2^2+\lVert\overline{\calC}_m^{(2)}\rVert^2_2 \bigr)^{1/2}+ 4\lVert \mathcal{E}(n,m)\rVert_2,
	\end{equation*}
	and whence
	\begin{align*}
	\lVert\overline\calC_{n+m}\rVert^2_2\le&\lVert \overline{\calC}_n^{(1)}\rVert_2^2+\lVert\overline{\calC}_m^{(2)}\rVert^2_2+8\bigl(\lVert \overline{\calC}_n^{(1)}\rVert_2^2+\lVert\overline{\calC}_m^{(2)}\rVert^2_2 \bigr)^{1/2}\lVert \mathcal{E}(n,m)\rVert_2\\&+16\lVert \mathcal{E}(n,m)\rVert^2_2.
	\end{align*}
	According to \cite[Lemmas 3.2 and  4.3]{CSS19},
	there is $c_4>0$ such that 
	\begin{equation}\label{E:h_d}\lVert \overline{\calC}_n\rVert_2\le c_4\sqrt{n}\quad\text{and}\quad \lVert \mathcal{E}(n,m)\rVert_2\le c_4h(n+m),\qquad n,m\ge1.
	\end{equation} 
	 Since the index of regular variation of $h(n)$ is strictly smaller than $1/2$, we arrive at
	\begin{align*}
	\lVert\overline\calC_{n+m}\rVert^2_2\le&\lVert \overline{\calC}_n^{(1)}\rVert_2^2+\lVert\overline{\calC}_m^{(2)}\rVert^2_2+8c_4^2\sqrt{n+m}h(n+m)+16c_4^2(h(n+m))^2\\
	&\le \lVert \overline{\calC}_n^{(1)}\rVert_2^2+\lVert\overline{\calC}_m^{(2)}\rVert^2_2+c_5\sqrt{n+m}h(n+m),
	\end{align*} for some $c_5>0$ large enough. 
	
	Analogously as above we have
	\begin{align*}\bigl(\lVert \overline{\calC}_n^{(1)}\rVert_2^2+\lVert\overline{\calC}_m^{(2)}\rVert^2_2 \bigr)^{1/2}
	&=
	\lVert \overline{\calC}_n^{(1)}+\overline{\calC}_m^{(2)}\rVert_2\\
	&\le \lVert\overline\calC_{n+m}\rVert_2+ \lVert \overline{\calC}_{n+m}-(\overline{\calC}_n^{(1)}+\overline{\calC}_m^{(2)})\rVert_2\\
	&\le \lVert\overline\calC_{n+m}\rVert_2+4\lVert \mathcal{E}(n,m)\rVert_2,
	\end{align*}
	which implies
	\begin{equation*}
	\lVert \overline{\calC}_n^{(1)}\rVert_2^2+\lVert\overline{\calC}_m^{(2)}\rVert^2_2\le \lVert\overline\calC_{n+m}\rVert^2_2+8\lVert\overline\calC_{n+m}\rVert_2\lVert \mathcal{E}(n,m)\rVert_2+16\lVert \mathcal{E}(n,m)\rVert_2^2.
	\end{equation*}
	Using \eqref{E:h_d} (and properties of $h(n)$) we conclude that 
	\begin{align*}
	\lVert \overline{\calC}_n^{(1)}\rVert_2^2+\lVert\overline{\calC}_m^{(2)}\rVert^2_2&\le \lVert\overline\calC_{n+m}\rVert^2_2+8c_4^2\sqrt{n+m} h(n+m)+16c_4^2(h(n+m))^2\\
	&\le \lVert\overline\calC_{n+m}\rVert^2_2+c_5\sqrt{n+m} h(n+m).
	\end{align*}
	We write
	\begin{equation*}
	x_n=\Var(\calC_n)=\lVert\overline{\calC}_n\rVert_2^2\quad\text{and}\quad b_n=c_5\sqrt{n} h(n),\qquad n\ge1. \end{equation*}
	We have shown that \begin{equation*}
	x_n+x_m-b_{n+m}\le x_{n+m}\le x_n+x_m+b_{n+m},\qquad n,m\ge1,
	\end{equation*}
	and from \cite[Lemmas 4.3 and 5.3]{CSS19} we know that 
	\begin{equation*}
	\lim_{n\nearrow\infty}\frac{x_n}{n}=\sigma^2>0.
	\end{equation*} 
	Take $n=m=2^{k-1}l$ for  $k,l\ge1$. Then one easily checks that
	\begin{equation*}
	\left|\frac{x_{2^kl}}{2^kl}-\frac{x_{2^{k-1}l}}{2^{k-1}l}\right|\le \frac{b_{2^kl}}{2^kl},\qquad k,l\ge1.
	\end{equation*}
	Next, observe that
	\begin{equation*}
	\sum_{k=1}^\infty\left(\frac{x_{2^kl}}{2^kl}-\frac{x_{2^{k-1}l}}{2^{k-1}l}\right)=\lim_{N\nearrow\infty}\sum_{k=1}^N\left(\frac{x_{2^kl}}{2^kl}-\frac{x_{2^{k-1}l}}{2^{k-1}l}\right)=\sigma^2-\frac{x_l}{l},\qquad l\ge1
	\end{equation*}
	and whence
	\begin{equation*}
	\left|\frac{x_n}{n}-\sigma^2\right|=\left|\sum_{k=1}^\infty\left(\frac{x_{2^kn}}{2^kn}-\frac{x_{2^{k-1}n}}{2^{k-1}n}\right)\right|\le\sum_{k=1}^\infty\frac{b_{2^kn}}{2^kn},\qquad n\ge1.
	\end{equation*}
	By recalling the definition of $\{b_n\}_{n\ge1}$, we conclude that
	\begin{equation*}
	\left|\frac{x_n}{n}-\sigma^2\right|\le\frac{c_5}{\sqrt{n}}\sum_{k=1}^\infty\frac{h(2^kn)}{2^{k/2}},\qquad n\ge1.
	\end{equation*}
	\textit{Case (i).} 
For $\Lambda\in(0,1/2)$ we apply \eqref{s_calC} and obtain 
	\begin{equation*}
	\APS{\frac{x_n}{n} - \sigma^2} \le \frac{c_5}{\sqrt{n}} \sum_{k = 1}^{\infty} \frac{(2^k n)^{1/2 - \Lambda} s(2^k n)}{2^{k/2}} = c_5 n^{-\Lambda} \sum_{k = 1}^{\infty} 2^{-k\Lambda} s(2^k n).
	\end{equation*}
	Since $s(n)$ is slowly varying, according to \cite[Theorem 1.5.6]{BGT_book}, there is a constant $c_6 > 0$ such that  $s(2^k n) \le c_6 2^{k\Lambda / 2} s(n)$ for all $k, n \ge 1$. Hence
	\begin{equation*}
	\APS{\frac{x_n}{n} - \sigma^2} \le c_5 c_6 n^{-\Lambda} s(n) \sum_{k = 1}^{\infty} 2^{-k\Lambda / 2} = c_7 n^{-\Lambda} s(n), \qquad n \ge 1,
	\end{equation*}
	which yields
	\begin{equation*}
	\aps{x_n - \sigma^2 n} \le c_7 n^{1- \Lambda} s(n)=c_7n^{1/2}h(n), \qquad n \ge 1.
	\end{equation*}
	\textit{Case (ii).} 
	If $\Lambda\ge1/2$, then $h(n)$ is slowly varying and thus 
	there is $c_8>0$ such that $h(2^kn)\le c_8 2^{k/4}h(n)$ for all $k,n\ge1$. This yields
	\begin{equation*}
	\left|\frac{x_n}{n}-\sigma^2\right|\le\frac{c_5c_8h(n)}{\sqrt{n}}\sum_{k=1}^\infty2^{-k/4}=c_9\frac{h(n)}{\sqrt{n}},\qquad n\ge1.
	\end{equation*}  Consequently
	\begin{equation*}
	|x_n-\sigma^2n|\le c_9n^{1/2}h(n),\qquad n\ge1,
	\end{equation*}
	and the proof is finished.
\end{proof}

\begin{remark}\label{RM:Var}
The following inequalities
\begin{align*}
\calC (0,n] \leq \calC_n \leq \calC (0,n] + \mathrm{Cap}(\{0\})
\end{align*}
entail 
$\barcalC_n -G^{-1}\leq \barcalC(0,n] \leq \barcalC_n + G^{-1}$,
where $G= G(0,0)=1/ \mathrm{Cap}(\{0\})$. Hence
\begin{align*}
\vert \barcalC (0,n] \vert \leq \vert\barcalC_n \vert + G^{-1}.
\end{align*}
Next, by taking squares of both sides of the above inequality and then expectations we obtain
\begin{align*}
\Var (\calC (0,n]) \leq \Var (\calC_n) + 2G^{-1}\bbE \big[ \vert \barcalC_n \vert \big] +  G^{-2}
\leq 
\Var (\calC_n) + 2G^{-1}\sqrt{\Var (\calC_n)} + G^{-2},
\end{align*}
where in the last inequality we applied Jensen's inequality. In view of Lemma \ref{LM:3.2} we infer that under assumptions A\ref{A1}, A\ref{A2} and $d/\alpha >5/2$ there exist constants  $c_1,c_2,c_3>0$ such that
\begin{align*}
\Var (\calC (0,n]) \leq \sigma^2n + c_1\sqrt{n}h(n) +2G^{-1}\big(\sigma^2n +c_2\sqrt{n}h(n)\big)^{1/2} + G^{-2}
\leq \sigma^2n + c_3\sqrt{n}h(n),
\end{align*}
where in the last inequality we used the fact that $h(n)$ is regularly varying of index less than $1/2$ if $d/\alpha >5/2$.
\end{remark}

\subsection*{The cardinality process}
We start with a decomposition of the range which goes back to Le Gall \cite{LeGall-French} and which was later applied in \cite[Corollary 2.1]{Asselah_Zd} to handle the capacity of the range.
For $a,b\in[0,\infty)$, $a\le b$, we use notation $\calR_a =\calR_{\floor{a}}$, $\calR [a,b]=\{S_a,\ldots , S_{b}\}$ and $\calR_{a-1}=\emptyset$ if $a<1$. It holds
\begin{align}\label{eq:123}
\begin{split}
|\calR_n|&=|\calR_{n/2}\cup\calR[n/2,n]|=|\calR_{n/2}\cup\calR[n/2,n]-S_{\lfloor n/2\rfloor}|\\&=|(\calR_{n/2}-S_{\lfloor n/2\rfloor})\cup(\calR[n/2,n]-S_{\lfloor n/2\rfloor})|=|\calR_{n/2}^{(1)}\cup\calR_{n/2}^{(2)}|\\&=|\calR_{n/2}^{(1)}|+|\calR_{n/2}^{(2)}|-|\calR_{n/2}^{(1)}\cap\calR_{n/2}^{(2)}|.
\end{split}
\end{align}
The Markov property implies that $\calR_{n/2}^{(1)}$ and $\calR_{n/2}^{(2)}$ are independent, and that $\calR_{n/2}^{(2)}$ is equal in law to $\calR_{\floor{n/2}}$ (or $\calR_{\floor{n/2+1}}$). By the symmetry of $\{S_n\}_{n\ge0}$ we have that $\calR_{n/2}^{(1)}$ has the same law as $\calR_{\floor{n/2}}$.  
By $I_n$ we denote the number of intersection points of two independent copies of our random walk up to time $n$, that is, $I_n=|\calR_{n}\cap\widetilde{\calR}_{n}|$
with $\widetilde{\calR}_{n}$ being independent of $\calR_{n}$, and having the same law.  
According to \cite[Remark on p.\ 666]{LeGall-Rosen} there is a constant $c_1>0$ such that
\begin{equation*}\label{abc}
\mathbb{E}[I_n]\le c_1 H(n),
\end{equation*}
where the function $H$ comes from \eqref{Function_H}.
%We fix the constant
%$\Delta = d/\alpha - 3/2$ and observe that if $\Delta \geq 1/2$ then $H(n)$ is slowly varying, see \cite[Lemma 2.2]{LeGall-Rosen},  and 
%if $\Delta\in (0, 1/2)$ we represent the function $H(n)$ in the form
%\begin{equation*}\label{eq:Fd_when_d_small}
%H(n) = n^{1/2 - \Delta} s(n),
%\end{equation*} 
%where $s(n)$ is a slowly varying function. In this case $H(n)$ is evidently regularly varying of index smaller than $1/2$.

By \cite[Theorem 4.4]{LeGall-Rosen},  assumption \textbf{A}\ref{A1} implies that if $d/\alpha > 3/2$ then the sequence $\{\Var(\calR_n)/n\}_{n \ge 1}$ converges to a strictly positive limit $\tilde{\sigma}^2$. We derive the second order term of this asymptotics.

\begin{lemma}\label{lm:var_asymp-card}
	Assume \textbf{A}\ref{A1} and $d/\alpha >3/2$. It holds that
	\begin{equation*}
	\Var(\calR_n) = 
	\tilde{\sigma}^2 n + O(\sqrt{n}\, H(n)).
	\end{equation*}
\end{lemma}
\begin{proof}
	Similarly as in \eqref{eq:123} we can easily show that for any $n,m\geq 1$,
	\begin{equation}\label{E:card_decomp}|\calR_n^{(1)}|+|\calR_m^{(2)}|-I_{n + m}\le |\calR_{n + m}|\le  |\calR_n^{(1)}|+|\calR_m^{(2)}|,
	\end{equation}
	where  $I_{n}=|\calR^{(1)}_{n}\cap\calR^{(2)}_{n}|$. 
	By taking expectation in \eqref{E:card_decomp} and then by subtracting those two relations, 
	\begin{equation*}
	\overline{\calR}^{(1)}_n + \overline{\calR}^{(2)}_m - I_{n + m} \le \overline{\calR}_{n + m} \le \overline{\calR}^{(1)}_n + \overline{\calR}^{(2)}_m + \bbE[I_{n + m}],
	\end{equation*}
	which implies
	\begin{align*}
	|\overline{\calR}_{n + m}-(\overline{\calR}_n^{(1)}+\overline{\calR}_m^{(2)})|
	\le I_{n + m}+\mathbb{E}[I_{n + m}]
	\le I_{n + m}+\lVert I_{n + m}\rVert_2.
	\end{align*}
	There is $c_4 > 0$ such that
	\begin{equation}\label{E:F_d}
	\lVert I_n\rVert_2\le c_4H(n)\quad\text{and}\quad 
	\lVert \overline{\calR}_n\rVert_2\le c_4\sqrt{n},\qquad n\ge1,
	\end{equation}
	see \cite[Remark after Corollary 3.2 and Theorem 4.4]{LeGall-Rosen}.
	Applying these estimates we can proceed similarly as in  \Cref{LM:3.2}. As $d > 3\alpha / 2$, $H(n) \le c_5\sqrt{n}$ for some $c_5>0$ and all $n\ge1$, and we obtain that for 
	\begin{equation*}
	x_n = \Var(\calR_n) = \norm{\overline{\calR}_n}_2^2 \qquad \text{and} \qquad b_n = c_6\sqrt{n} H(n),\qquad n\ge1, \end{equation*}
	for some $c_6>0$, we have
	\begin{equation*}
	x_n+x_m-b_{n+m}\le x_{n+m}\le x_n+x_m+b_{n+m},\qquad n,m\ge1,
	\end{equation*}
	and from \cite[Theorem 4.4]{LeGall-Rosen} we know that
	\begin{equation*}
	\lim_{n\nearrow\infty}\frac{x_n}{n}=\tilde{\sigma}^2>0.
	\end{equation*}
	Setting $n=m=2^{k-1}l$ for  $k,l\ge1$, we arrive at
	\begin{equation*}
	\left|\frac{x_n}{n}-\tilde{\sigma}^2\right|\le\frac{c_6}{\sqrt{n}}\sum_{k=1}^\infty\frac{H(2^kn)}{2^{k/2}},\qquad n\ge1.
	\end{equation*}
 Again, by a similar reasoning as in \Cref{LM:3.2}, the result follows.
\end{proof}

\subsection*{Acknowledgement}
We thank the referees for their helpful comments and suggestions.
N.\ Sandri\'c was supported by  
\textit{Alexander-von-Humboldt Foundation} under project No.\ HRV 1151902 HFST-E, and \textit{Croatian Science Foundation} under project 8958.
S.\ \v Sebek was supported by 
\textit{Austrian Scienc Fund FWF} under project P31889-N35, and \textit{Croatian Science Foundation} under project 4197.

\bibliographystyle{abbrv}
\bibliography{LIL}

\end{document}